\documentclass[11pt]{article}
\usepackage[english]{babel}
\usepackage{amsmath,amsthm,amssymb}
\usepackage{mathrsfs}
\usepackage{bbm}
\usepackage{amsfonts}
\usepackage[margin=0.8in]{geometry}
\usepackage{graphicx}
\usepackage{caption}
\usepackage{subcaption}
\usepackage{wasysym}
\usepackage{enumerate}
\usepackage{algorithm2e}
\usepackage{tabularx}
\usepackage{color}
\usepackage{wrapfig}
\usepackage{todonotes}

\newtheorem{thm}{Theorem}[section]
\newtheorem{ass}[thm]{Assumption}
\newtheorem{cor}[thm]{Corollary}
\newtheorem{lem}[thm]{Lemma}
\newtheorem{prop}[thm]{Proposition}
\newtheorem*{hyp*}{Hypothesis}
\theoremstyle{definition}
\newtheorem{defn}[thm]{Definition}

\theoremstyle{rem}
\newtheorem{rem}[thm]{Remark}

\numberwithin{equation}{section}

\newcommand{\R}{\mathbb R}
\newcommand{\eps}{\varepsilon}

\newcommand{\bbD}{\mathbb D}
\newcommand{\bbF}{\mathbb F}

\newcommand{\bbS}{\mathbb S}

\newcommand{\bbT}{\mathbb T}

\newcommand{\mcA}{\mathcal{A}}

\newcommand{\mcB}{\mathcal{B}}
\newcommand{\mcC}{\mathcal C}
\newcommand{\mcD}{\mathcal D}
\newcommand{\mcE}{\mathcal E}

\newcommand{\mcF}{\mathcal F}
\newcommand{\mcI}{\mathcal I}

\newcommand{\mcJ}{\mathcal J}

\newcommand{\mcT}{\mathcal T}
\newcommand{\mcP}{\mathcal P}

\newcommand{\mcU}{\mathcal U}

\newcommand{\mcY}{\mathcal Y}
\newcommand{\PrM}{\mathfrak{P}}

\newcommand{\E}{\mathbb{E}}
\newcommand{\Prob}{\mathbb{P}}

\newcommand{\vecv}{\bold{v}}
\newcommand{\ett}{\mathbbm{1}}
\newcommand{\cadlag}{c\`adl\`ag~}
\newcommand{\cadlagSTOP}{c\`adl\`ag}

\DeclareMathOperator*{\argmin}{arg\,min}

\newcommand{\ie}{\textit{i.e.\ }}
\newcommand{\eg}{\textit{e.g.\ }}
\newcommand{\etal}{\textit{et.~al.\ }}

\begin{document}

\title{Non-Markovian Impulse Control Under Nonlinear Expectation\footnote{This work was supported by the Swedish Energy Agency through grant number 48405-1}}

\author{Magnus Perninge\footnote{M.\ Perninge is with the Department of Physics and Electrical Engineering, Linnaeus University, V\"axj\"o,
Sweden. e-mail: magnus.perninge@lnu.se.}} %
\maketitle
\begin{abstract}
We consider a general type of non-Markovian impulse control problems under adverse non-linear expectation or, more specifically, the zero-sum game problem where the adversary player decides the probability measure. We show that the upper and lower value functions satisfy a dynamic programming principle (DPP). We first prove the dynamic programming principle (DPP) for a truncated version of the upper value function in a straightforward manner. Relying on a uniform convergence argument then enables us to show the DPP for the general setting. Following this, we use an approximation based on a combination of truncation and discretization to show that the upper and lower value functions coincide, thus establishing that the game has a value and that the DPP holds for the lower value function as well. Finally, we show that the DPP admits a unique solution and give conditions under which a saddle-point for the game exists.

As an example, we consider a stochastic differential game (SDG) of impulse versus classical control of path-dependent stochastic differential equations (SDEs).
\end{abstract}

\section{Introduction}
We solve a robust impulse control problem where the aim is to find an impulse control, $u^*$, that solves
\begin{align}
\inf_{u\in\mcU}\mcE^u\Big[\varphi((\tau_i,\beta_i)_{i=1}^N) + \sum_{j=1}^N c((\tau_i,\beta_i)_{i=1}^j)\Big],\label{ekv:minmax-form}
\end{align}
where $\mcE^u[\cdot]:=\sup_{\Prob\in\mcP(u)}\E^\Prob[\cdot]$ is a non-linear expectation and the terminal reward $\varphi$ is a random function that maps impulse controls $u=(\tau_j,\beta_j)_{j=1}^\infty$ ($N:=\max\{j:\tau_j< T\}$) to values of the real line and is measurable with respect to $\mcF\otimes\mcB(\bar\mcD)$, where $\mcB(\bar\mcD)$ is the Borel $\sigma$-field of the space $\bar\mcD$ where the control $u$ takes values and $\mcF$ is the $\sigma$-field generated by the canonical process on the space of continuous trajectories starting at 0. The intervention cost $c$ is defined similarly to $\varphi$ but is in addition assumed to be progressively measurable in the last intervention-time and bounded from below by a positive constant.

In fact, we take this formulation one step further by showing that, under mild conditions, the zero-sum game where we play an impulse control while the adversary player (nature) chooses a probability measure has a value when allowing our adversary to play strategies, \ie that
\begin{align}
\inf_{u\in\mcU}\sup_{\Prob\in\mcP(u)}\E^\Prob\Big[\varphi((\tau_i,\beta_i)_{i=1}^N) + \sum_{j=1}^N c((\tau_i,\beta_i)_{i=1}^j)\Big] = \sup_{\Prob^S\in\mcP^S}\inf_{u\in\mcU}\E^{\Prob^S(u)}\Big[\varphi((\tau_i,\beta_i)_{i=1}^N) + \sum_{j=1}^N c((\tau_i,\beta_i)_{i=1}^j)\Big]\label{ekv:game-form}
\end{align}
where $\mcP^S$ is the set of non-anticipative strategies mapping impulse controls to probability measures, and derive additional conditions under which a saddle-point, $(u^*,\Prob^{*,S})\in \mcU\times\mcP^S$, for the game exists.

Our approach relies on the tower property of non-linear expectations discovered in \cite{NutzHandel} and applied in \cite{NutzZhang15} to solve an optimal stopping problem under adverse non-linear expectation. In this regard we need to assume that $\varphi$ and $c$ are uniformly bounded and uniformly continuous under a suitable metric.

To indicate the applicability of the results we consider the special case when \eqref{ekv:minmax-form} corresponds to the stochastic differential game (SDG) of impulse versus classical control
\begin{align}
\inf_{u\in\mcU}\sup_{\alpha\in \mcA}\E\Big[\int_0^T\phi(s,X^{u,\alpha}_s)ds+\psi(X^{u,\alpha}_T) - \sum_{j=1}^N\ell(\tau_j,X^{(\tau_i,\beta_i)_{i=1}^{j-1},\alpha}_{\tau_j},\beta_j)\Big], \label{ekv:sdg-form}
\end{align}
where $\mcA$ is a set of classical controls and $X^{u,\alpha}$ solves an impulsively-continuously controlled path-dependent stochastic differential equation (SDE) that implements $u$ in feedback form. To assure sufficient regularity in this setting, we impose an additional $L^2$-Lipschitz condition on the coefficients of the SDE. In particular, we will see that this SDG corresponds to setting
\begin{align}\label{ekv:phiSFDE}
\varphi(\omega,u)&=\int_0^T\phi(s,(\mcI^u(\omega))_s)ds+\psi((\mcI^u(\omega))_T),
\end{align}
and
\begin{align}\label{ekv:cSFDE}
c(\omega,(\tau_i,\beta_i)_{i=1}^j)&=\ell(\tau_j(\omega),(\mcI^{(\tau_i,\beta_i)_{i=1}^{j-1}}(\omega))_{\tau_j(\omega)},\beta_j(\omega)),
\end{align}
where for each $u\in\mcU$, $\mcI^u$ is a map from the space of continuous paths to the space of \cadlag paths that adds the impulses in $u$ to the, otherwise, continuous trajectory.

The main contributions of the present work are threefold. First, we show that the game \eqref{ekv:game-form} has a solution when the sets $(\mcP(u):u\in\mcU)$ satisfy standard conditions translated to our setting and the functions $\varphi$ and $c$ are bounded and uniformly continuous. Second, we extract a saddle-point under additional weak-compactness assumptions on the family $(\mcP(u):u\in\mcU)$. Finally, we give a set of conditions under which the cost/reward pair defined by \eqref{ekv:phiSFDE}-\eqref{ekv:cSFDE} satisfies the assumptions in the first part of the paper, enabling us to show that the path-dependent SDG of classical versus impulse control in \eqref{ekv:sdg-form} has a value.\\


\textbf{Related literature} The optimal stopping problem under adverse non-linear expectation was considered by Nutz and Zhang \cite{NutzZhang15} and by Bayraktar and Yao in \cite{BayraktarYao14}, where the latter allows for a slightly more general setting not having to assume a uniform bound on the rewards. Nevertheless, as explained above, \cite{NutzZhang15} is based on the tower property of non-linear expectations developed in \cite{NutzHandel} and is, therefore, more closely related to the present work.

Non-Markovian impulse control under standard (linear) expectation was first considered by Hamad\`ene \etal in~\cite{DjehiceImpulse}, where it was assumed that the impulses do not affect the dynamics of the underlying process. This approach was extended to incorporate delivery lag in~\cite{Hdhiri} and, more recently, also to an infinite horizon setting in~\cite{DjehicheInfHorImp}. A different approach to non-Markovian impulse control was initiated in \cite{SwitchElephant} and then further developed in \cite{JonteSFDE} where interconnected Snell envelopes indexed by controls were used to find solutions to problems with impulsively controlled path-dependent SDEs. We mention also the general formulation of impulse controls in~\cite{ImpulsIH}, which can be seen as a linear expectation version of the present work, and the work on impulse control of path-dependent SDEs under $g$-expectation (see~\cite{Peng04}) and related systems of backward SDEs (BSDEs) in \cite{SeqBSDE}. Although the latter work considers a path-dependent SDG of impulse versus classical control, the classical control only enters the drift term. Effectively this corresponding to the situation when the set $\mcP(u)$ in our framework is dominated and the extension to non-dominated sets would have to go through the incorporation of second order BSDEs (2BSDEs) (see \cite{2bsde07,STZ2bsde}).

The idea of having one player implement a strategy in the zero-sum game setting was first proposed by Elliot and Kalton~\cite{ElliotKalton72} to counter the unrealistic idea that one of the players have to give up their control to the opponent in games of control versus control. The approach was combined with the theory of viscosity solutions to find a representation of the upper and lower value functions in deterministic differential games as solutions to Hamilton-Jacobi-Bellman-Isaacs (HJBI) equations by Evans and Souganidis \cite{Evans84}. Using a discrete time approximation technique, this was later translated to the stochastic setting by Flemming and Souganidis \cite{FlemSoug89}. Notable is that, while the framework of \cite{ElliotKalton72}, that has been the prevailing formulation in the literature since its introduction, assumes that the first player to act always implements a strategy, our formulation only allows the adversary player to implement a strategy when she acts first, while our impulse control is always open loop control. The reason that our approach is still successful lies in the weak formulation of the game, effectively turning the impulse control into a feedback control when we turn to the SDG in the latter part of the paper. Moreover, our approach avoids the asymmetric information structure that results from implementing the game formulation in~\cite{ElliotKalton72} which enables use to derive saddle-points.

Related to the SDG of impulse versus classical control that we consider is the work of Azimzadeh~\cite{Azimzadeh19} when the intervention costs are deterministic and by Bayraktar \etal \cite{BayraktarRobust} that considers a robust impulse control problem when the impulse control is of switching type. Both of these are restricted to the Markovian setting. The latter implements the switching control in feedback form while the classical control is open loop. In this sense it is probably the closes work to the present one that can be found in the literature. However, the approach used in \cite{BayraktarRobust} is based on the Stochastic Perron Method of Bayraktar and S\^irbu~\cite{BayraktarSirbu12,BayraktarSirbu13} (see also \cite{Sirbu14} for another application to SDGs) that do not easily translate to a path-dependent setting. The setup in \cite{Azimzadeh19} was later extended in \cite{ImpVClass} to allow for stochastic intervention costs and $g$-nonlinear expectation. However, this extension also falls within the Markovian framework and uses standard BSDEs to define the cost/reward.

Previous works on non-Markovian SDGs in the standard framework of classical versus classical control can be found in Pham and Zhang \cite{PhamZhang14} and Possama\"i \etal \cite{Possamai20}, where the latter uses a weak-formulation of feedback versus feedback control and shows that the game has a value under the Isaac's condition and uniqueness of solutions to the corresponding path-dependent HJBI-equation~\cite{EkrenPathDep_p1,EkrenPathDep_p2}. We remark that an interesting further development of the present work would be to consider the corresponding path-dependent quasi-variational inequalities, reminiscent of the relation between Markovian impulse control problems and classical quasi-variational inequalities~\cite{BensLionsImpulse}.\\

\textbf{Outline} In the next section we give some preliminary definitions, recall some prior results (in particular, the solution to the optimal stopping problem under non-linear expectation in \cite{NutzZhang15}) and give an immediate extension of the tower property in~\cite{NutzHandel} to our setting. Then, in Section~\ref{sec:DynPP} we show that a type of dynamic programming principle holds for the upper value function. In Section~\ref{sec:value}, we use an approximation routine applied to both value functions to conclude that the game has a value, \ie that \eqref{ekv:game-form} holds. In Section~\ref{sec:verif-thm} it is shown that our DPP admits a unique solution and conditions are given under which an optimal pair can be extracted from the value function. Finally, in Section~\ref{sec:sdg} we relate the result to path-dependent SDGs of impulse versus classical control.

\section{Preliminaries\label{sec:prel}}
\subsection{Notation}
Throughout, we shall use the following notation, where we set a fixed horizon $T\in (0,\infty)$:
\begin{itemize}
  \item We fix a positive integer $d$, define the sample space $\Omega:=\{\omega\in C(\R_+\to\R^d):\omega_0=0\}$ and set $\Lambda:=[0,T]\times\Omega$. For $t\in [0,T]$ we introduce the pseudo-norm $\|\omega\|_t:=\sup_{s\in[0,t]}|\omega_s|$ and extend the corresponding distance to $\Lambda$ by defining
      \begin{align}
        \mathbf d_{\Lambda}[(t,\omega),(t',\omega')]:=|t'-t|+\|\omega'_{\cdot\wedge t'}-\omega_{\cdot\wedge t}\|_T.
      \end{align}
  \item For $(t,\omega)\in\Lambda$, we let $\Omega^{t,\omega}:=\{\omega'\in\Omega:\omega'|_{[0,t]}=\omega|_{[0,t]}\}$.
  \item The set of all probability measures on $\Omega$ equipped with the topology
of weak convergence, \ie the weak topology induced by the bounded continuous functions on
$\Omega$, is denoted $\PrM(\Omega)$.
  \item We let $B$ denote the canonical process, \ie $B_t(\omega):=\omega_t$ and denote by $\Prob_0$ the probability measure under which $B$ is a Brownian motion and let $\E$ be expectation with respect to $\Prob_0$.
  \item $\bbF:=\{\mcF_t\}_{0\leq t\leq T}$ is the natural (raw) filtration generated by $B$ and $\bbF^*:=\{\mcF^*_t\}_{0\leq t\leq T}$, where $\mcF^*_t$ is the universal completion of $\mcF_t$. 
  \item We let $\mcT$ be the set of all $\bbF$-stopping times and for each $\eta\in\mcT$ we let $\mcT_\eta$ be the set of all $\tau\in\mcT$ such that $\tau(\omega)\geq\eta(\omega)$ for all $\omega\in\Omega$. For fixed $t\in [0,T]$, we let $\mcT^t$ be all $\tau\in\mcT_t$ such that $\omega\mapsto\tau(\omega)$ is independent of $\omega|_{[0,t]}$.
  \item For $\kappa\geq0$, we let $D^\kappa:=\{(t_j,b_j)_{j=1}^{\kappa}: 0\leq t_1\leq\cdots\leq t_\kappa\leq T,\, b_j\in U\}$, where $U$ is a compact subset of $\R^d$, and set $\mcD:=\cup_{\kappa\geq 0}D^\kappa$ and $\bar\mcD:=\mcD\cup \{(t_j,b_j)_{j= 1}^\infty: 0\leq t_1\leq\cdots\leq T,\, b_j\in U\}$. Moreover, for $0\leq t\leq s\leq T$, we let $\mcD_{[t,s]}$ (resp. $\bar\mcD_{[t,s]}$) be the subset of $(t_j,b_j)_{j=1}^\kappa\in\mcD$ (resp. $\bar\mcD$) with $t_j\subset [t,s]$ for $j=1,\ldots,\kappa$.
  \item We let $\mcU$ be the set of all $u=(\tau_j,\beta_j)_{j=1}^N$, where $(\tau_j)_{j=1}^\infty$ is a non-decreasing sequence of $\bbF$-stopping times, $\beta_j$ is a $\mcF_{\tau_j}$-measurable r.v.~taking values in $U$ and $N:=\max\{j:\tau_j< T\}$.  Moreover, for $k\geq 0$ we let $\mcU^k:=\{u\in\mcU:\,N\leq k\}$. We interchangeably refer to elements $u\in \mcU$ (resp. $u\in\mcU^k$) as $u=(\tau_j,\beta_j)_{j=1}^N$ and $u=(\tau_j,\beta_j)_{j=1}^\infty$ (resp. $u=(\tau_j,\beta_j)_{j=1}^k$).
  \item For $t\in[0,T]$, we let $\mcU_t$ (resp $\mcU^k_t$) be the subset of $\mcU$ (resp $\mcU^k$) with all controls for which $\tau_1\geq t$ and denote by $\mcU^{t}$ the subset of $\mcU_t$ with all controls $u$ such that $\omega\mapsto u(\omega)$ is independent of $\omega|_{[0,t]}$.
  \item For $u\in\mcU\cup(\cup_{j=1}^\infty\mcU^j)$ and $k\geq 0$, we let $[u]_k:=(\tau_j,\beta_j)_{j=1}^{N\wedge k}$. Furthermore, for $t\geq 0$ we let $N(t):=\max\{j\geq 0:\tau_j\leq t\}$ and define $u_t:=[u]_{N(t)}$ and $u^t:=(\tau_j,\beta_j)_{j= N(t)+1}^N$.
  \item We denote by $\emptyset\in \mcU$ the control with no interventions, \ie $N=0$ implying that $\tau_1=T$.
  \item For each $\kappa\geq 0$, we introduce the pseudo-distance $\mathbf d_k$ on $\Lambda\times D^\kappa$ as
  \begin{align*}
    \mathbf d_\kappa[(t,\omega,\vecv),(t',\omega',\vecv')]:=\|\omega'_{\cdot\wedge t'}-\omega_{\cdot\wedge t}\|_T+\sum_{j=1}^\kappa\|\omega'_{\cdot\wedge t'_j\wedge t'}-\omega_{\cdot\wedge t_j\wedge t}\|_T+|\vecv'-\vecv|.
  \end{align*}
\end{itemize}

We stretch the definition of uniform continuity to maps on $\Lambda\times\mcD$ as follows:
\begin{defn}
We say that a map $\mcY:\Lambda\times\mcD\to\R$ is uniformly continuous if for each $\kappa\geq 0$, the map $\mcY:\Lambda\times D^\kappa\to\R$ is uniformly continuous under the distance $\mathbf d_\kappa$.
\end{defn}

We define the concatenation of $\vecv=(t_j,b_j)_{j=1}^\kappa\in D^\kappa$ and $\vecv'=(t'_j,b'_j)_{j\geq 1}\in\bar\mcD$ as
\begin{align*}
  \vecv\circ\vecv':=((t_1,b_1),\ldots,(t_\kappa,b_\kappa),(t'_1\vee t_\kappa,b'_1),\ldots,(t'_j\vee t_\kappa,b'_j),\ldots).
\end{align*}
For technical reasons we extend the composition to pairs $(\vecv,\vecv')$ when $\vecv$ has infinite length by letting $\vecv\circ\vecv'=\vecv$ in this case. The extension allows us to decompose any control $u\in\mcU$ as $u=u_\tau\circ u^\tau$.

For $\tau\in\mcT$ and $\omega,\omega'\in\Omega$ we introduce the composition on $\Omega$ as
\begin{align*}
  (\omega \otimes_\tau\omega')_s:=\omega_s\ett_{[0,\tau(\omega))}(s)+(\omega_{\tau(\omega)}+\omega'_{s-\tau(\omega)})\ett_{[\tau(\omega),T]}(s),\qquad \forall s\in[0,T]
\end{align*}
and for $f:\Omega\to\R$ we set $f^{\tau,\omega}(\omega'):=f(\omega \otimes_\tau\omega')$.

The results we present rely on the notion of regular conditional probability: Any $\Prob\in\PrM(\Omega)$ has a regular conditional probability distribution $(\Prob^\omega_\tau)_{\omega\in\Omega}$ given $\mcF_\tau$ satisfying
\begin{align*}
  \Prob^\omega_\tau[\{\omega'\in\Omega:\omega'=\omega\text{ on }[0,\tau(\omega)]\}]=1
\end{align*}
for all $\omega\in\Omega$ (see \eg page 34 in \cite{StroockVaradhan}). We define the probability measure $\Prob^{\tau,\omega}\in\PrM(\Omega)$ as
\begin{align*}
  \Prob^{\tau,\omega}[A]:=\Prob^\omega_\tau[\omega\otimes_\tau A],
\end{align*}
for all $A\in\mcF$, where $\omega\otimes_\tau A:=\{\omega\otimes_\tau\omega':\omega'\in A\}$.


\begin{defn}
We let $\{\mcP(t,\omega,\vecv\circ u):(t,\omega,\vecv,u)\in \cup_{t\in [0,T]}\{t\}\times\Omega\times \bar\mcD_{[0,t]} \times \mcU_t\}$ be a family of subsets of $\PrM(\Omega)$ such that
\begin{align*}
  \mcP(t,\omega,\vecv\circ u)=\mcP(t,\omega',\vecv\circ u)\quad \text{when }\omega|_{[0,t]}=\omega'|_{[0,t]}.
\end{align*}
Moreover, we let $\mcP(\tau,\omega,u):=\mcP(\tau(\omega),\omega,u_{\tau(\omega)}(\omega)\circ u^{\tau(\omega)})$ and set $\mcP(u):=\mcP(0,\omega,u)$.
\end{defn}

We recall that a subset of a Polish space is analytic if it is the image of a Borel subset of another Polish space under a Borel map and that an $\bar R$-valued function $f$ is upper semi-analytic if the set $\{f> c\}$ is analytic for each $c\in \R$.

For $\Prob\in\PrM(\Omega)$ we let $\E^\Prob_t[f](\omega):=\E^\Prob[f^{t,\omega}]$ (here $\E^\Prob$ is expectation with respect to $\Prob$) and define the non-linear expectation
\begin{align*}
  \mcE^u_\tau[f](\omega):=\sup_{\Prob\in\mcP(\tau,\omega,u)}\E^\Prob_\tau[f](\omega)
\end{align*}
for all $(\tau,\omega,u)\in\mcT\times\Omega\times \mcU$ and all upper semi-analytic functions $f:\Omega\to\bar\R$. Then, $\mcE^u_\tau[f]$ is $\mcF_\tau^*$-measurable and upper semi-analytic~\cite{NutzHandel}.

The idea is that the adversary player, given a trajectory $\omega|_{[0,t]}$ and a sequence of impulses $\vecv\in\mcD_{[0,t]}$ chooses a probability measure on $\Omega^{t,\omega}$ to maximize \eqref{ekv:minmax-form}. In this regard we introduce the set of non-anticipative maps from controls to probability measures:
\begin{defn}\label{def:strat}
We denote by $\mcP^S$ the set of non-anticipative maps $\Prob^S:\mcU\to\PrM(\Omega)$ mapping $u\in\mcU$ to $\Prob=\Prob^S(u)\in\mcP(u)$. By non-anticipativity, we mean that if $u_{\tau-}=\tilde u_{\tau-}$ for some $\tau\in\mcT$ and $u,\tilde u\in\mcU$, then $\Prob^S(u)[A]=\Prob^S(\tilde u)[A]$ for all $A\in\mcF_\tau$.

Moreover, for $(t,\omega)\in\Lambda$ we let $\mcP^S(t,\omega)$ denote the set of all non-anticipative maps $\Prob^S:\mcU\to\PrM(\Omega)$ such that $\Prob^S(u)\in\mcP(t,\omega,u)$. Often, we suppress dependence of $\omega$ and write $\mcP^S_t$ for $\mcP^S(t,\omega)$.
\end{defn}

\subsection{Assumptions}
\begin{ass}\label{ass:mcP}
For each $u\in\mcU$ we assume that $\mcP(u)$ is non-empty. Moreover, $\Prob\in\mcP(u)$ if for each $k\geq 0$, there is a $\Prob'\in \mcP([u]_k)$ with $\Prob'=\Prob$ on $\mcF_{\tau_k}$.
\end{ass}

\begin{ass}\label{ass:mcE-general}
For $(t,\omega)\in\Lambda$, $\vecv\in \bar\mcD_{[0,t]}$ and $u\in\mcU_t$, we assume that $\mcP(t,\omega,\vecv\circ u)$ satisfies for any stopping time $\tau\in\mcT_t$ and $\Prob\in\mcP(t,\omega,\vecv\circ u)$ (with $\theta:=\tau^{t,\omega}-t$):
\begin{enumerate}[a)]
  \item\label{ass:mcE-gen-graph} The graph $\{(\Prob',\omega'):\omega'\in\Omega,\,\Prob'\in\mcP(\tau,\omega',\vecv\circ u)\}\subset\PrM(\Omega)\times\Omega$ is analytic.
  \item We have $\Prob^{\theta,\omega'}\in\mcP(\tau,\omega\otimes_{t}\omega',\vecv\circ u)$ for $\Prob$-a.e.~$\omega'$.
  \item\label{ass:mcE-gen-Kernel} If $\nu:\Omega\to\PrM(\Omega)$ is an $\mcF_{\theta}$-measurable kernel and $\nu(\omega')\in\mcP(\tau,\omega\otimes_{t}\omega',\vecv\circ u)$ for $\Prob$-a.e.~$\omega'$, then the measure defined by
  \begin{align}\label{ekv:kernel-def}
    \Prob'[A]:=\int\!\!\!\int (\ett_A)^{\theta,\tilde \omega}(\omega')\nu(d\omega';\tilde\omega)\Prob[d\tilde\omega]
  \end{align}
  for any $A\in\mcF$, belongs to $\mcP(t,\omega,\vecv\circ u)$.
  \item\label{ass:mcE-gen-consist} Let $\tilde u\in\mcU_t$ be such that $u_{\tau-}=\tilde u_{\tau-}$, then there is a $\Prob'\in\mcP(t,\omega,\vecv\circ \tilde u)$ such that $\Prob'=\Prob$ on $\mcF_{\theta}$.
\end{enumerate}
\end{ass}
In the above assumption, conditions \ref{ass:mcE-gen-graph})-\ref{ass:mcE-gen-Kernel}) are standard (see \eg \cite{NeufeldNutz13}) and basically means that dynamic programming holds.  The last condition implies that the non-anticipativity postulated in Definition~\ref{def:strat} is achievable. 

\begin{ass}\label{ass:costs}
The functions $\varphi:\Omega\times \mcD\to\R$ and $c:\Omega\times \mcD\to [\delta,\infty)$ (with $\delta>0$) are uniformly bounded, \ie there is a constant $C_0>0$ such that
\begin{align*}
|\varphi(\omega,\vecv)|+|c(\omega,\vecv)|\leq C_0,
\end{align*}
for all $(\omega,\vecv)\in\Omega\times\bar\mcD$. For each $\kappa\geq 0$, there is a modulus of continuity function $\rho_{\varphi,\kappa}$ such that
\begin{align*}
|\varphi(\omega',\vecv')-\varphi(\omega,\vecv)|\leq \rho_{\varphi,\kappa}(\mathbf d_\kappa[(T,\omega,\vecv),(T,\omega',\vecv')])
\end{align*}
for all $\omega,\omega'\in\Omega$ and $\vecv,\vecv'\in D^\kappa$ with $t_j\leq t'_j$. Moreover, $c(u)$ is $\mcF_{\tau_N}$-measurable and for each $\kappa\geq 0$ there is a modulus of continuity function $\rho_{c,\kappa}$ such that
\begin{align*}
|c(\omega',\vecv')-c(\omega,\vecv)|\leq \rho_{c,\kappa}(\mathbf d_\kappa[(t_\kappa,\omega,\vecv),(t'_\kappa,\omega',\vecv')]),
\end{align*}
where $t_\kappa$ and $t'_\kappa$ are the times of the last interventions in $\vecv$ and $\vecv'$, respectively. 
Finally, for any $\vecv\in\mcD$ and $b\in U$ we have
\begin{align}\label{ekv:no-imp@end}
  \varphi(\omega,\vecv)< \varphi(\omega,\vecv\circ (T,b))+c(\omega,\vecv\circ (T,b))
\end{align}
for all $\omega\in\Omega$.
\end{ass}


\begin{rem}\label{rem:@end}
In view of \eqref{ekv:no-imp@end} it is never optimal to intervene on the system at time $T$. In light of this and to simplify notation later on we will assume that $\varphi$ is such that $\varphi(\omega,\vecv\circ(T,b))=\varphi(\omega,\vecv)$ for all $\omega\in\Omega$, $\vecv\in\mcD$ and $b\in U$.
\end{rem}

\begin{ass}\label{ass:mcE-consist}
For each $i\geq 0$, there is a modulus of continuity $\rho_{\mcE,i}$ such that for all $\kappa,k\geq 0$, $\vecv',\vecv\in D^\kappa$,  $t\in[0,T]$ and $\omega, \omega'\in\Omega $ and $u\in\mcU^k_t$,  there is a $u_\omega\in\mcU^k_{t}$ such that
\begin{align}\nonumber
&\mcE^{\vecv'\circ u}_{t}\Big[\varphi(\vecv'\circ u) + \sum_{j=1}^N c(\vecv'\circ [u]_j)\Big](\omega')-\mcE^{\vecv\circ u_\omega}_t\Big[\varphi(\vecv\circ u_\omega) + \sum_{j=1}^{N_\omega} c(\vecv\circ [u_\omega]_j)\Big](\omega)
\\
&\leq \rho_{\mcE,\kappa+k}(\mathbf d_\kappa[(t,\omega,\vecv),(t,\omega',\vecv')])\label{ekv:mcE-cont}
\end{align}
and $(\omega,\omega')\mapsto u_\omega(\omega')$ is $\mcF_t\otimes\mcF$-measurable.
\end{ass}

\begin{ass}\label{ass:mcE-bound} For each modulus of continuity $\rho$ there is a modulus of continuity $\rho'$ such that for any $\tau\in\mcT$, we have
$\sup_{\Prob\in\mcP(t,\omega,u)}\E^\Prob(\rho(\eps+\sup_{s\in [\tau,\tau+\eps]}|B_s|))\leq \rho'(\eps)$ for all $(t,\omega,u)\in\Lambda\times\mcU$.
\end{ass}



\subsection{The tower property for our family of nonlinear expectations}
The following trivial extension of Theorem 2.3 in \cite{NutzHandel} follows immediately from the above assumptions:

\begin{prop}\label{prop:tower}
For each $\tau,\eta\in\mcT$, with $\tau\leq\eta$, $\vecv\in\mcD$ and $u\in\mcU_\tau$ and $\tilde u\in\mcU_\eta$ we have
\begin{align*}
  \mcE^{\vecv\circ u}_\tau\big[\mcE^{\vecv\circ u\circ \tilde u}_\eta\big[f\big]\big]=\mcE^{\vecv\circ u\circ \tilde u}_\tau\big[f\big]
\end{align*}
for each upper semi-analytic function $f:\Omega\to\bar\R$.
\end{prop}

\noindent \emph{Proof.} We first note that since $\mcE^{\vecv\circ u\circ \tilde u}$ satisfies Assumption 2.1 in  \cite{NutzHandel}, Theorem 2.3 of the same article gives that
\begin{align*}
  \mcE^{\vecv\circ u \circ \tilde u}_\tau\big[\mcE^{\vecv\circ u\circ \tilde u}_\eta\big[f\big]\big]=\mcE^{\vecv\circ u\circ \tilde u}_\tau\big[f\big]
\end{align*}
Moreover, since $(\vecv\circ u \circ \tilde u)_{\eta-}=(\vecv\circ u)_{\eta-}$, Assumption~\ref{ass:mcE-general}.(\ref{ass:mcE-gen-consist}) gives that for each $\Prob\in \mcP(t,\omega,\vecv\circ u \circ \tilde u)$, there is a $\Prob'\in \mcP(t,\omega,\vecv\circ u)$ such that $\Prob=\Prob'$ on $\mcF_\eta$ (and thus on $\mcF^*_\eta$), and vice versa. Finally, as $\mcE^{\vecv\circ u\circ \tilde u}_\eta\big[f\big]$ is $\mcF^*_\eta$-measurable the result follows.\qed\\

\subsection{Optimal stopping under non-linear expectation}
We recall the following result on optimal stopping under non-linear expectation $\mcP,\,\mcE$ (satisfying the assumptions on, say, $\mcP(\emptyset),\,\mcE^{\emptyset}$ above).

\begin{thm}\label{thm:NutzZhang}{(Nutz and Zhang (2015)\cite{NutzZhang15})}
Assume that the process $(X_t)_{0\leq t\leq T}$ has \cadlag paths, is progressively measurable and bounded and satisfies
\begin{align*}
  X_{t'}(\omega')-X_{t}(\omega)\leq \rho_X(\mathbf d_{\Lambda}[(t,\omega),(t',\omega')]),
\end{align*}
for some modulus of continuity function $\rho_X$ and all $(t,\omega),(t',\omega')\in\Lambda$ with $ t\leq t'$. Then, $Y_t:=\inf_{\tau\in\mcT_t}\mcE_t[X_\tau]$ satisfies:
\begin{enumerate}[i)]
  \item Let $\tau^*:=\inf\{s\geq t:Y_s=X_s\}$, then $\tau^*\in\mcT_t$ and $Y_t:=\mcE_t[X_{\tau^*}]$. Moreover, $Y_{\cdot\wedge\tau^*}$ is a $\Prob$-supermartingale for any $\Prob\in\mcP$.
  \item The game has a value: $\inf_{\tau\in\mcT}\mcE[X_\tau]=\sup_{\Prob\in\mcP}\inf_{\tau\in\mcT}\E^\Prob[X_\tau]$.
  \item If $\mcP(t,\omega)$ is weakly compact for each $(t,\omega)\in\Lambda$, then there is a $\Prob^*\in\mcP$ such that  $Y_0=\E^{\Prob^*}[X_{\tau^*}]=\inf_{\tau\in\mcT} \E^{\Prob^*}[X_{\tau}]$.
\end{enumerate}
\end{thm}

\section{A dynamic programming principle\label{sec:DynPP}}
We define the upper value function as\footnote{Clearly, we also have $Y^{\vecv}_t(\omega)=\inf_{u\in\mcU^t}\mcE^{\vecv\circ u}\big[\varphi(\vecv\circ u) + \sum_{j=1}^N c(\vecv\circ [u]_j)\big](\omega)$.}
\begin{align}
Y^{\vecv}_t(\omega):=\inf_{u\in\mcU_t}\mcE^{\vecv\circ u}\big[\varphi(\vecv\circ u) + \sum_{j=1}^N c(\vecv\circ [u]_j)\big](\omega)
\end{align}
for all $\omega\in\Omega$ and set $Y=Y^\emptyset$. Moreover, we define the lower value function as
\begin{align}
Z^{\vecv}_t(\omega):=\sup_{\Prob\in\mcP^S_t}\inf_{u\in\mcU}\E_t^{\Prob(\vecv\circ u)}\big[\varphi(\vecv\circ u) + \sum_{j=1}^N c(\vecv\circ [u]_j)\big](\omega)
\end{align}
for all $\omega\in\Omega$. 

\begin{rem}
As noted in the introduction, it may appear as though the setup is somewhat asymmetric as the optimization problem for the lower value function contains a strategy, whereas the corresponding problem for the upper value function does not. However, as the impulse controls in $\mcU$ are $\bbF$-adapted and the opponent controls the probability measure (effectively deciding the likelihoods for different trajectories), this can be seen as the impulse player implementing a type of non-anticipative strategy as well. This becomes more evident when turning to the application in Section \ref{sec:sdg}.
\end{rem}

In this section we will concentrate on the upper value function and the main result is the following dynamic programming principle:
\begin{thm}\label{thm:dynP}
The map $Y$ is bounded, uniformly continuous and satisfies the recursion
\begin{align}\label{ekv:dynp-Y}
    Y^{\vecv}_t(\omega)=\inf_{\tau\in\mcT^t}\mcE^{\vecv}_t\Big[\ett_{[\tau=T]}\varphi(\vecv) + \ett_{[\tau<T]}\inf_{b\in U}\{Y^{\vecv\circ(\tau,b)}_\tau+c(\vecv\circ(\tau,b))\}\Big](\omega)
  \end{align}
for all $(t,\omega,\vecv)\in \Lambda\times\mcD$.
\end{thm}

\begin{rem}
As a consequence, the map $(t,\omega,\vecv)\mapsto Y^{\vecv}_t(\omega)$ is Borel measurable.
\end{rem}


The proof of Theorem~\ref{thm:dynP} is given through a sequence of lemmata where the main obstacle that we need to overcome is to show uniform continuity of the map $(t,\omega,\vecv)\mapsto Y^{\vecv}_t(\omega)$. This will be obtained through a uniform convergence argument and we introduce the truncated upper value function defined as
\begin{align}
Y^{\vecv,k}_ t(\omega):=\inf_{u\in\mcU^k_t}\mcE^{\vecv\circ u}\Big[\varphi(\vecv\circ u) + \sum_{j=1}^N c(\vecv\circ [u]_j)\Big](\omega)
\end{align}
for all $(t,\omega,\vecv)\in\Lambda\times\mcD$ and $k\geq 0$. Note that in the definition of $Y^{\vecv,k}$, the impulse controller is restricted to using a maximum of $k$ impulses. The following approximation result is central:

\begin{lem}\label{lem:unif-conv}
There is a $C>0$ such that
\begin{align*}
  Y^{\vecv}_t(\omega)\geq Y^{\vecv,k}_t(\omega)-C/k
\end{align*}
for all $(t,\omega,\vecv)\in \Lambda\times\mcD$ and each  $k\geq 0$.
\end{lem}

\noindent \emph{Proof.} We first note that
\begin{align*}
  \inf_{u\in\mcU_t}\mcE^{\vecv\circ u}_t\Big[\varphi(\vecv\circ u) +\sum_{j=1}^N c(\vecv\circ [u]_{j})\Big]\leq \mcE^{\vecv}_t\Big[\varphi(\vecv) \Big]\leq C_0.
\end{align*}
Now, as
\begin{align*}
  \mcE^{\vecv\circ u}_t\Big[\varphi(\vecv\circ u) +\sum_{j=1}^N c(\vecv\circ [u]_{j})\Big]\geq-C_0 + \mcE^{\vecv\circ u}_t\big[N \delta\big]
\end{align*}
for any $u\in\mcU_t$, this implies that in the infimum of \eqref{ekv:dynp-Y} we can restrict our attention to impulse controls for which
\begin{align}\label{ekv:u-sens}
  k\mcE^{\vecv\circ u}_t\big[\ett_{[N\geq k]}\big]\leq \mcE^{\vecv\circ u}_t\big[N\big]\leq 2C_0/\delta.
\end{align}
for each $k\geq 0$, since any impulse control not satisfying \eqref{ekv:u-sens} is dominated by $u=\emptyset$. For any $u\in\mcU_t$ satisfying \eqref{ekv:u-sens} we have,
\begin{align*}
   &\mcE^{\vecv\circ [u]_k}_t\Big[\varphi(\vecv\circ [u]_k) +\sum_{j=1}^{N\wedge k} c(\vecv\circ [u]_{j})\Big]
   \\
   &=\mcE^{\vecv\circ [u]_k}_t\Big[(\ett_{[\tau_{k+1} = T]}+\ett_{[\tau_{k+1}<T]})\varphi(\vecv\circ [u]_k) +\sum_{j=1}^{N\wedge k} c(\vecv\circ [u]_{j})\Big]
   \\
   &\leq \mcE^{\vecv\circ [u]_k}_t\Big[\ett_{[\tau_{k+1} = T]}\varphi(\vecv\circ u) +\sum_{j=1}^{N\wedge k} c(\vecv\circ [u]_{j})\Big]+\mcE^{\vecv\circ [u]_k}_t\Big[\ett_{[\tau_{k+1}<T]}\varphi(\vecv\circ [u]_k)\Big]
   \\
   &\leq \mcE^{\vecv\circ [u]_k}_t\Big[\ett_{[\tau_{k+1} = T]}\varphi(\vecv\circ u) +\sum_{j=1}^{N\wedge k} c(\vecv\circ [u]_{j})\Big]+2C_0^2/(k+1)\delta.
\end{align*}
To arrive at the last inequality above, we have used that since $([u]_k)_{\tau_{k+1}-}=u_{\tau_{k+1}-}$, Assumption~\ref{ass:mcE-general}.\ref{ass:mcE-gen-consist}) gives that
\begin{align*}
  \mcE^{\vecv\circ [u]_k}_t\big[\ett_{[\tau_{k+1}<T]}\varphi(\vecv\circ [u]_k)\big]&\leq C_0 \mcE^{\vecv\circ [u]_k}_t\big[\ett_{[\tau_{k+1}<T]}\big]
  \\
  &=C_0\mcE^{\vecv\circ u}_t\big[\ett_{[\tau_{k+1}<T]}\big]
   \\
  &=C_0\mcE^{\vecv\circ u}_t\big[\ett_{[N\geq k+1]}\big]
\end{align*}
from which the inequality is immediate by \eqref{ekv:u-sens}. Now, as $\ett_{[\tau_{k+1} = T]}\varphi(\vecv\circ u) +\sum_{j=1}^{N\wedge k} c(\vecv\circ [u]_{j})$ is $\mcF_{\tau_{k+1}}$-measurable we can again use Assumption~\ref{ass:mcE-general}.\ref{ass:mcE-gen-consist}) to find that
\begin{align*}
    \mcE^{\vecv\circ [u]_k}_t\Big[\ett_{[\tau_{k+1} = T]}\varphi(\vecv\circ u) +\sum_{j=1}^{N\wedge k} c(\vecv\circ [u]_{j})\Big]= \mcE^{\vecv\circ u}_t\Big[\ett_{[\tau_{k+1} = T]}\varphi(\vecv\circ u) +\sum_{j=1}^{N\wedge k} c(\vecv\circ [u]_{j})\Big].
\end{align*}
Consequently,
\begin{align*}
   &\mcE^{\vecv\circ [u]_k}_t\Big[\varphi(\vecv\circ [u]_k) +\sum_{j=1}^{N\wedge k} c(\vecv\circ [u]_{j})\Big]-\mcE^{\vecv\circ u}_t\Big[\varphi(\vecv\circ u) +\sum_{j=1}^N c(\vecv\circ [u]_{j})\Big]
   \\
   &\leq \mcE^{\vecv\circ u}_t\Big[(\ett_{[\tau_{k+1} = T]}-1)\varphi(\vecv\circ u) -\sum_{j=k+1}^{N} c(\vecv\circ [u]_{j})\Big]+2C_0^2/(k+1)\delta
   \\
   &\leq C_0\mcE^{\vecv\circ u}_t\Big[\ett_{[\tau_{k+1} < T]}\Big]+2C_0^2/(k+1)\delta
   \\
   &\leq 4C_0^2/(k+1)\delta.
\end{align*}
Since $u$ was an arbitrary impulse control satisfying \eqref{ekv:u-sens} and $[u]_k\in\mcU^k_t$, the result follows.\qed\\

\begin{lem}\label{lem:Yk-prop}
For each $k,\kappa\geq 0$ and $t\in[0,T]$, the map $(\omega,\vecv)\mapsto Y^{\vecv,k}_t(\omega):\Omega\times D^\kappa\to\R$ is uniformly continuous with a modulus of continuity that is independent of $t$.
\end{lem}

\noindent \emph{Proof.} We have
\begin{align*}
Y^{\vecv',k}_{t}(\omega')-Y^{\vecv,k}_t(\omega)&\leq\sup_{u\in\mcU^k_{t}}\Big\{\mcE^{\vecv'\circ u}_t\big[\varphi(\vecv'\circ u) + \sum_{j=1}^N c(\vecv'\circ [u]_j)\big](\omega')\\
&\quad-\mcE^{\vecv\circ u_{\omega}}_t\big[\varphi(\vecv\circ u_{\omega}) + \sum_{j=1}^N c(\vecv\circ [u_{\omega}]_j)\big](\omega)\Big\}
\end{align*}
with $u_\omega$ as in Assumption \ref{ass:mcE-consist}. This immediately gives that
\begin{align*}
Y^{\vecv',k}_{t}(\omega')-Y^{\vecv,k}_t(\omega)\leq \rho_{\mcE,\kappa+k}(\textbf d_\kappa[(t,\omega,\vecv),(t,\omega',\vecv')])
\end{align*}
by the same assumption.\qed\\

As usual, the dynamic programming principle for $Y^{\cdot,k}$ will be proved by leveraging regularity and we introduce the following partition of the set $\Lambda\times U$:
\begin{defn}\label{def:discretization}
For $\eps>0$:
\begin{itemize}
  \item We let $n\geq 0$ be the smallest integer such that $2^{-n}T\leq\eps$, set $n_t^\eps:=2^{n}+1$ and introduce the discrete set $\bbT^\eps:=\{t^\eps_j:t^\eps_j=(j-1)2^{-n}T,j=1,\ldots,n^\eps_t\}$. Moreover, we define $\bar t^\eps_{i}:=t^\eps_{i}$ for $i=1,\ldots,n_t^\eps-1$ and $\bar t^\eps_{n_t^\eps}:=T+1$.
  \item We then let $(E^\eps_{i,j})_{1\leq i\leq n^\eps_t}^{1\leq j}$ be such that $(E^\eps_{i,j})_{j\geq 1}$ forms a partition of $\Omega$ with $E^\eps_{i,j}\in\mcF_{t^\eps_i}$ and $\|\omega-\omega'\|_{t^\eps_i}\leq \eps$ for all $\omega,\omega'\in E^\eps_{i,j}$. We let $(\omega^\eps_{i,j})_{1\leq i\leq n^\eps_t}^{j\leq 1}$ be a sequence with $\omega^\eps_{i,j}\in E^\eps_{i,j}$.
  \item Finally, we let $(U^\eps_{l})_{1\leq l\leq n^\eps_U}$ be a Borel-partition of $U$ such that the diameter of $U^\eps_l$ does not exceed $\eps$ for $l=1,\ldots,n^\eps_U$ and let $(b^\eps_l)_{l=1}^{n^\eps_U}$ be a sequence with $b^\eps_l\in U^\eps_l$ and denote by $\bar U^\eps:=\{b^\eps_1,\ldots,b^\eps_{n^\eps_U}\}$ the corresponding dicretization of $U$.
\end{itemize}

\end{defn}

In the next lemma we show that for any $\eps>0$, the impulse size can be chosen $\eps$-optimally.

\begin{lem}\label{lem:eps-attaind}
Let $g:\Lambda\times\mcD\to\R$ be bounded and uniformly continuous. Then, for each $k\geq 0$, $\eps>0$, $u\in\mcU^k$ and $\tau\in \mcT_{\tau_k}$, there is a $\mcF_{\tau}$-measurable random variable $\beta$ with values in $U$, such that
\begin{align*}
    \sup_{\Prob\in\mcP(t,\omega,v)}\E^\Prob\big[g(\tau,u\circ(\tau,\beta))-\inf_{b\in U}g(\tau,u\circ(\tau,b))\big]\leq \eps
  \end{align*}
for all $(t,\omega,v)\in\Lambda\times \mcU$.
\end{lem}

\noindent \emph{Proof.} \textbf{Step 1.} We first prove the result for $u=\vecv\in D^k$. For arbitrary $\eps_1>0$ there is, by properties of the supremum, a double sequence $(b_{i,j})_{1\leq i\leq n^{\eps_1}_t}^{j\geq 1}$ with $b_{i,j}\in U$ such that
\begin{align*}
  g(t^{\eps_1}_i,\omega^{\eps_1}_{i,j},\vecv\circ(t^{\eps_1}_i,b_{i,j}))\leq \inf_{b\in U}g(t^{\eps_1}_i,\omega^{\eps_1}_{i,j},\vecv\circ(t^{\eps_1}_i,b))+\eps/2
\end{align*}
for $i=1,\ldots,n^{\eps_1}_t$ and all $j\geq 1$. Consequently, letting $\beta^\vecv(\omega'):=\sum_{i=1}^{n_t-1}\sum_{j\geq 0}b_{i,j}\ett_{[t^\eps_i,\bar t^\eps_{i+1})}(\tau)\ett_{E^{\eps_1}_{i,j}}(\omega')$ we find that
\begin{align*}
  g(\tau(\omega'),\omega',\vecv\circ(\tau(\omega'),\beta^\vecv(\omega')))&\leq \inf_{b\in U}g(\tau(\omega'),\omega',\vecv\circ(\tau(\omega'),b))
  \\
  &\quad+2\rho_{g}(3\eps_1+\|\omega'_{(\tau(\omega')+\cdot)\wedge T}-\omega'_{\tau(\omega')}\|_{\eps_1})+\eps/2,
\end{align*}
where $\rho_g$ is the modulus of continuity of $g$ and the result, for deterministic $u$, follows from Assumption~\ref{ass:mcE-bound}.\\

\textbf{Step 2.} We turn to the general setting with $u\in\mcU^k$ and let
\begin{align*}
    \hat u:=\Big(\sum_{i=1}^{n^{\eps_1}_t-1}t^{\eps_1}_{i+1}\ett_{[t^{\eps_1}_{i},\bar t^{\eps_1}_{i+1})}(\tau_j),\sum_{i=1}^{n^{\eps_1}_U} b^{\eps_1}_i\ett_{U^{\eps_1}_i}(\beta_j)\Big)_{j=1}^k.
 \end{align*}
By definition, it follows that
\begin{align*}
   |g(s,\omega',\hat u(\omega')\circ  (s,b))-g(s,\omega',u(\omega')\circ(s,b))|\leq \rho_{g}((2k+3)\eps_1+2\sum_{j=1}^k\|\omega'_{\tau_j(\omega')+\cdot}-\omega'_{\tau_j(\omega')}\|_{\eps_1})
\end{align*}
for all $(s,\omega',b)\in\Lambda\times U$. Let $\beta:=\beta^{\hat u}$. Then $\beta$ is $\mcF_\tau$-measurable and we find that
\begin{align*}
    &g(\tau(\omega'),\omega',u(\omega')\circ(\tau(\omega'),\beta(\omega')))-\inf_{b\in U}g(\tau(\omega'),\omega',u(\omega')\circ(\tau(\omega'),b))
    \\
    &\leq g(\tau(\omega'),\omega',\hat u(\omega')\circ(\tau(\omega'),\beta(\omega')))-\inf_{b\in U}g(\tau(\omega'),\omega',\hat u(\omega')\circ(\tau(\omega'),b))
    \\
    &\quad+2\rho_{g}((2k+3)\eps_1+2\sum_{j=1}^k\|\omega'_{\tau_j(\omega')+\cdot}-\omega'_{\tau_j(\omega')}\|_{\eps_1})
    \\
    &\leq 2\rho_{g}(3\eps_1+\|\omega'_{(\tau(\omega')+\cdot)\wedge T}-\omega_{\tau(\omega)}\|_{\eps_1})+ 2\rho_{g}((2k+3)\eps_1+2\sum_{j=1}^k\|\omega_{\tau_j(\omega')+\cdot}-\omega_{\tau_j(\omega')}\|_{\eps_1})+\eps/2.
  \end{align*}
Finally, by choosing a concave modulus of continuity that dominates $\rho_g$ the result follows from using Assumption~\ref{ass:mcE-bound}, taking nonlinear expectations on both sides and choosing $\eps_1$ sufficiently small.\qed\\

\begin{lem}\label{lem:dynp-Yk}
For each $k\geq 0$, the map $(t,\omega,\vecv)\mapsto Y^{\vecv,k}_t(\omega)$ is bounded and uniformly continuous. Moreover, the family $(Y^{\cdot,k})_{k\geq 0}$ satisfies the recursion
\begin{align}\label{ekv:dynp-Yk}
    Y^{\vecv,k}_t(\omega)=\inf_{\tau\in\mcT^t}\mcE^{\vecv}_t\Big[\ett_{[\tau=T]}\varphi(\vecv) + \ett_{[\tau<T]}\inf_{b\in U}\{Y^{\vecv\circ(\tau,b),k-1}_\tau+c(\vecv\circ(\tau,b))\}\Big](\omega),
  \end{align}
for each $k\geq 1$ and $(t,\omega,\vecv)\in\Lambda\times\mcD$.
\end{lem}

\noindent \emph{Proof.} We note that $|Y^{\vecv,k}_t|\leq C_0$ and so boundedness clearly holds. The remaining assertions will follow by induction in $k$ and we note by a simple argument that $Y^{\cdot,0}$ is uniformly continuous and bounded. In fact, for $0\leq t\leq s\leq T$ we have by the tower property that
\begin{align}\nonumber
  |Y^{\vecv,0}_t(\omega)-Y^{\vecv,0}_s(\omega)|&=|\mcE^\vecv_t\big[Y^{\vecv,0}_s\big](\omega)-Y^{\vecv,0}_s(\omega)|
  \\
  &\leq \mcE^\vecv_t\big[|Y^{\vecv,0}_s-Y^{\vecv,0}_s(\omega)|\big](\omega)\nonumber
  \\
  &\leq \sup_{\Prob\in\mcP(t,\omega,\vecv)}\E^\Prob\big[\rho_{\mcE,\kappa}(\|(\omega\otimes_t B)_{\cdot}-\omega_\cdot\|_{s})\big]\nonumber
  \\
  &\leq \sup_{\Prob\in\mcP(t,\omega,\vecv)}\E^\Prob\big[\rho_{\mcE,\kappa}(\mathbf d_{\Lambda}[(s,\omega),(t,\omega)]+\|B\|_{s-t})\big]\nonumber
  \\
  &\leq \rho'(\mathbf d_{\Lambda}[(s,\omega),(t,\omega)])\label{ekv:Y0-cont}
\end{align}
that together with Lemma~\ref{lem:Yk-prop} yields the desired continuity.

We thus assume that for each $\kappa\geq 0$, the map $Y^{\cdot,k-1}:\Lambda\times D^\kappa\to\R$ is uniformly continuous with modulus of continuity $\rho_{\kappa,k-1}$ and prove that then \eqref{ekv:dynp-Yk} holds and that this, in turn, implies uniform continuity of $Y^{\cdot,k}$. Let us consider the right-hand side of \eqref{ekv:dynp-Yk} that we denote $\hat Y^{\vecv,k}_t(\omega)$. By our induction assumption and Theorem \ref{thm:NutzZhang}.\emph{i)} we have,
\begin{align*}
    \hat Y^{\vecv,k}_t(\omega)=\mcE^{\vecv}_t\Big[\ett_{[\tau_1^*=T]}\varphi(\vecv) + \ett_{[\tau_1^*<T]}\inf_{b\in U}\{Y^{\vecv\circ(\tau_1^*,b),k-1}_{\tau^*_1}+c(\vecv\circ(\tau_1^*,b))\}\Big](\omega),
  \end{align*}
where $\tau_1^*=\inf\big\{s\geq t: \hat Y^{\vecv,k}_s(\omega)=\inf_{b\in U}\{Y^{\vecv\circ(s,b),k-1}_{s}+c(\vecv\circ(s,b))\}\big\}\wedge T$.\\

\textbf{Step 1.} We show that $\hat Y\geq Y$. For this we fix $\kappa\geq 0$ and $\vecv\in D^\kappa$ and note that for each $\eps>0$, our induction assumption together with Lemma~\ref{lem:eps-attaind} implies the existence of a $\mcF_{\tau_1^*}$-measurable $\bar\beta^\eps_1$ with values in $U$ such that
\begin{align*}
    \hat Y^{\vecv,k}_t(\omega)&\geq \mcE^{\vecv}_t\Big[\ett_{[\tau_1^*=T]}\varphi(\vecv) + \ett_{[\tau_1^*<T]}\{Y^{\vecv\circ(\tau_1^*,\bar\beta^\eps_1),k-1}_{\tau_1^*}+c(\vecv\circ (\tau_1^*,\bar\beta^\eps_1))\}\Big](\omega) - \eps/2
    \\
    &\geq \mcE^{\vecv}_t\Big[\ett_{[\tau_1^*=T]}\varphi(\vecv) + \sum_{i=2}^{n_t^\eps}\ett_{[t^\eps_{i-1},t^\eps_{i})}(\tau_1^*) \sum_{l=1}^{n^\eps_U}\ett_{U^\eps_l}(\bar\beta^\eps_1)\{Y^{\vecv_{i,l},k-1}_{t^\eps_{i}}+c(\vecv_{i,l})
    \\
    &\quad-(\rho_{\kappa+1,k-1}+\rho_{c,\kappa+k})(\textbf d[(\tau_1^*,\cdot,\vecv\circ (\tau_1^*,\beta^\eps_1)),(t^\eps_{i},\cdot,\vecv_{i,l})])\}\Big](\omega) - \eps/2
    \\
    &\geq \mcE^{\vecv}_t\Big[\ett_{[\tau_1^*=T]}\varphi(\vecv) + \sum_{i=2}^{n_t^\eps}\ett_{[t^\eps_{i-1},t^\eps_{i})}(\tau_1^*) \sum_{l=1}^{n^\eps_U}\ett_{U^\eps_l}(\bar\beta^\eps_1)\{Y^{\vecv_{i,l},k-1}_{t^\eps_{i}}+c(\vecv_{i,l})\}\Big](\omega) - \rho'(\eps) - \eps/2,
\end{align*}
where $\vecv_{i,l}:=\vecv\circ (t^\eps_{i},b^\eps_l)$ and $\rho'$ is a modulus of continuity function. On the other hand, for each $(i,j)$, there is a $ u^\eps_{i,j,l}\in\mcU_{t^\eps_i}^{k-1}$ such that
\begin{align*}
  Y^{\vecv_{i,l},k-1}_{t^\eps_i}(\omega^\eps_{i,j})\geq \mcE^{\vecv_{i,l}\circ u^\eps_{i,j,l}}_{t^\eps_i}\Big[\varphi(\vecv_{i,l}\circ u^\eps_{i,j,l}) +\sum_{m=1}^{N^\eps} c(\vecv_{i,l}\circ [u^\eps_{i,j,l}]_{m})\Big](\omega^\eps_{i,j})-\eps/2.
\end{align*}
Moreover, Assumption~\ref{ass:mcE-consist} implies the existence of a $\mcF_{t^\eps_i}\times\mcF$-measurable map $(\omega,\omega')\mapsto u^{i,j,l}_{\omega}(\omega')$ with $u^{i,j,l}_{\omega}\in\mcU^{k-1}_{t^\eps_i}$ for all $\omega\in\Omega$, such that
\begin{align*}
  &\mcE^{\vecv_{i,l}\circ u^{i,j,l}_{\omega}}_{t^\eps_i}\Big[\varphi(\vecv_{i,l}\circ u^{i,j,l}_{\omega}) +\sum_{m=1}^{N^{i,j,l}_{\omega}} c(\vecv_{i,l}\circ [u^{i,j,l}_{\omega}]_{m})\Big](\omega)
  \\
  &\leq \mcE^{\vecv_{i,l}\circ u_{i,j,l}}_{t^\eps_i}\Big[\varphi(\vecv_{i,l}\circ u_{i,j,l}) +\sum_{m=1}^{N_{i,j,l}} c(\vecv_{i,l}\circ [u_{i,j,l}]_{m})\Big](\omega^\eps_{i,j}) + \rho_{\mcE,\kappa+k}(\eps).
\end{align*}
for all $\omega\in E^\eps_{i,j}$. By once again using uniform continuity of $Y^{\cdot,k-1}$ we have that
\begin{align*}
  |Y^{\vecv_{i,l},k-1}_{t^\eps_i}(\omega^\eps_{i,j})-Y^{\vecv_{i,l},k-1}_{t^\eps_i}(\omega)|\leq \rho_{\kappa+1,k-1}(\eps)
\end{align*}
for all $\omega\in E^\eps_{i,j}$, and we conclude that
\begin{align*}
  Y^{\vecv_{i,l},k-1}_{t^\eps_i}(\omega)\geq \mcE^{\vecv_{i,l}\circ u^{i,j,l}_{\omega}}_{t^\eps_i}\Big[\varphi(\vecv_{i,l}\circ u^{i,j,l}_{\omega}) +\sum_{m=1}^{N^{i,j,l}_{\omega}} c(\vecv_{i,l}\circ [u^{i,j,l}_{\omega}]_{m})\Big](\omega)- \rho_{\mcE,\kappa+k}(\eps)- \rho_{\kappa+1,k-1}(\eps) - \eps/2
\end{align*}
whenever $\omega\in E^\eps_{i,j}$. Letting $\tilde u:=\sum_{i= 2}^{n^\eps_t}\ett_{[t^\eps_{i-1},t^\eps_i)}(\tau^\eps_1)\sum_{j\geq 1}\ett_{E^\eps_{i,j}}(\omega)\sum_{l=1}^{n^\eps_U}\ett_{U^\eps_{l}}(\beta_1^\eps)u^{i,j,l}_{\omega}$, we find that
\begin{align*}
   Y^{\vecv_{i,l},k-1}_{t^\eps_{i}}(\omega)&\geq \mcE^{\vecv_{i,l}\circ \tilde u}_{t^\eps_{i}}\Big[\varphi(\vecv_{i,l}\circ \tilde u) +\sum_{j=1}^{\tilde N} c(\vecv_{i,l}\circ [\tilde u]_{j})\Big](\omega)- \rho_{\mcE,\kappa+k}(\eps)- \rho_{\kappa+1,k-1}(\eps)- \eps/2
\end{align*}
for all $\omega\in E^\eps_{i,j}$. We can combine the above impulse controls into
\begin{align*}
  u^\eps:=\Big(\sum_{i=2}^{n_t^\eps}t^\eps_{i}\ett_{[t^\eps_{i-1},t^\eps_{i})}(\tau_1^*)+T\ett_{[\tau^*_1=T]}, \sum_{l=1}^{n^\eps_U}b_l\ett_{U^\eps_l}(\bar\beta^\eps_1)\Big)\circ \tilde u\in\mcU^k_t.
\end{align*}
The tower property now gives that
\begin{align*}
\hat Y^{\vecv,k}_ t(\omega)\geq\mcE^{\vecv\circ u^\eps}_t\Big[\varphi(\vecv\circ u^\eps) +\sum_{j=1}^{N} c(\vecv\circ [u^\eps]_{j})\Big](\omega)-\eps-\rho(\eps)
\end{align*}
for some modulus of continuity $\rho$ and it follows that $\hat Y^{\vecv,k}_ t(\omega)\geq Y^{\vecv,k}_ t(\omega)$ since $\eps>0$ was arbitrary.\\

\textbf{Step 2.} We now show that $\hat Y\leq Y$. Pick $\hat u=(\hat \tau_j,\hat \beta_j)_{j=1}^k\in\mcU^k_t$ and note that
\begin{align*}
  \hat Y^{\vecv,k}_t(\omega)\leq \mcE^{\vecv}_t\Big[\ett_{[\hat \tau_1=T]}\varphi(\vecv) + \ett_{[\hat \tau_1<T]}\{Y^{\vecv\circ(\hat \tau_1,\hat \beta_1),k-1}_{\tau_1}+c(\vecv\circ(\hat \tau_1,\hat \beta_1))\}\Big](\omega).
\end{align*}
Moreover, for each $\omega\in\Omega$ we have
\begin{align*}
  Y^{\vecv\circ(\hat \tau_1,\hat \beta_1),k-1}_{\hat \tau_1}(\omega)\leq \mcE^{\vecv\circ \hat u}_{\hat \tau_1}\Big[\varphi(\vecv\circ \hat u) + \sum_{j=2}^{\hat N} c(\vecv \circ [\hat u]_{j})\Big](\omega).
\end{align*}
We conclude by the tower property that
\begin{align*}
  \hat Y^{\vecv,k}_t(\omega)&\leq \mcE^{\vecv}_t\Big[\ett_{[\hat \tau_1=T]}\varphi(\vecv) + \ett_{[\hat \tau_1<T]}\{\mcE^{\vecv\circ \hat u}_{\hat \tau_1}\Big[\varphi(\vecv\circ \hat u) + \sum_{j=2}^{\hat N} c(\vecv \circ [\hat u]_{j})\Big]+c(\vecv\circ(\hat \tau_1,\hat \beta_1))\}\Big](\omega)
  \\
  &= \mcE^{\vecv\circ \hat u}_{t}\Big[\varphi(\vecv\circ \hat u) + \sum_{j=1}^{\hat N} c(\vecv \circ [\hat u]_{j})\Big]
\end{align*}
and it follows that $\hat Y^{\vecv,k}_ t(\omega)\leq Y^{\vecv,k}_ t(\omega)$ since this time $\hat u\in \mcU^k_t$ was arbitrary.\\

\textbf{Step 3.} It remains to show that $Y^{\cdot,k}$ is uniformly continuous. As uniform continuity in $(\omega,\vecv)$ follows from Lemma~\ref{lem:Yk-prop} we only need to consider the time variable and find a modulus of continuity that is independent of $\vecv$. Let $0\leq t\leq s\leq T$ and note that by the preceding steps and (4.4) of \cite{NutzZhang15}
\begin{align*}
  |Y^{\vecv,k}_t(\omega)-Y^{\vecv,k}_s(\omega)|&\leq \mcE^\vecv_t\big[|Y^{\vecv,k}_s-Y^{\vecv,k}_s(\omega)|\big](\omega)
  \\
  &+\sup_{\tau\in\mcT^t}\mcE^\vecv_t\big[(Y^{\vecv,k}_s - \ett_{[\tau<T]}\inf_{b\in U}\{Y^{\vecv\circ(\tau,b),k-1}_{\tau}+c(\vecv\circ(\tau,b))\})\ett_{[\tau<s]}\big](\omega)
\end{align*}
By arguing as in \eqref{ekv:Y0-cont} we have
\begin{align*}
\mcE^\vecv_t\big[|Y^{\vecv,k}_s-Y^{\vecv,k}_s(\omega)|\big](\omega)&\leq \sup_{\Prob\in\mcP(t,\omega,\vecv)}\E^\Prob\big[\rho_{\mcE,\kappa+k}(\mathbf d_{\Lambda}[(s,\omega),(t,\omega)]+\|B\|_{s-t})\big]
\\
&\leq \rho(\mathbf d_{\Lambda}[(s,\omega),(t,\omega)])
\end{align*}
for some $\rho$ independent of $\vecv$. Concerning the second term we have, since
\begin{align*}
  Y^{\vecv,k}_s\leq \inf_{b\in U}\{Y^{\vecv\circ(s,b),k-1}_{s}+c(\vecv\circ(s,b))\}
\end{align*}
that
\begin{align*}
  &\mcE^\vecv_t\big[(Y^{\vecv,k}_s - \ett_{[\tau<T]}\inf_{b\in U}\{Y^{\vecv\circ(\tau,b),k-1}_{\tau}+c(\vecv\circ(\tau,b))\})\ett_{[\tau<s]}\big](\omega)
  \\
  &\leq \mcE^\vecv_t\big[(\inf_{b\in U}\{Y^{\vecv\circ(s,b),k-1}_{s}+c(\vecv\circ(s,b))\} - \inf_{b\in U}\{Y^{\vecv\circ(\tau,b),k-1}_{\tau}+c(\vecv\circ(\tau,b))\})\ett_{[\tau<s]}\big](\omega)
  \\
  &\leq \mcE^\vecv_t\big[\sup_{b\in U}\{Y^{\vecv\circ(s,b),k-1}_{s}-Y^{\vecv\circ(\tau,b),k-1}_{\tau}+c(\vecv\circ(s,b)) - c(\vecv\circ(\tau,b)) \}\ett_{[\tau<s]}\big](\omega)
  \\
  &\leq \sup_{\Prob\in\mcP(s,\omega,\vecv)}\E^\Prob\big[\rho_{\kappa+1,k-1}(\mathbf d_{\Lambda}[(s,\omega\otimes_s B),(t,\omega\otimes_s B)]+\sup_{\Prob\in\mcP(s,\omega,\vecv)}\E^\Prob\big[\rho_{c,\kappa+k}(\mathbf d_{\Lambda}[(s,\omega\otimes_s B),(t,\omega\otimes_s B)])
  \\
  &\leq \rho'(\mathbf d_{\Lambda}[(s,\omega),(t,\omega)])
\end{align*}
for some modulus of continuity $\rho'$ and we conclude that $(t,\omega)\mapsto Y^{\vecv,k}_t(\omega)$ is uniformly continuous with a modulus of continuity that does not depend on $\vecv$. Uniform continuity of the map $(t,\omega,\vecv)\mapsto Y^{\vecv,k}_t(\omega)$ then follows by Lemma~\ref{lem:Yk-prop}. This completes the induction-step.\qed\\

\noindent \emph{Proof of Theorem~\ref{thm:dynP}.} The sequence $(Y^{\vecv,k}_{t}(\omega))_{k\geq 0}$ is bounded, uniformly in $(t,\omega,\vecv)$, and since $\mcU^k\subset\mcU^{k+1}$ it follows by Lemma~\ref{lem:unif-conv} that $Y^{\vecv,k}_{t}(\omega)\searrow Y^{\vecv}_{t}(\omega)$. Moreover, again by Lemma~\ref{lem:unif-conv} the convergence is uniform and Lemma~\ref{lem:dynp-Yk} gives that $Y^{\vecv}_{t}(\omega)$ is bounded and uniformly continuous. Concerning the dynamic programming principle \eqref{ekv:dynp-Y}, we have by \eqref{ekv:dynp-Yk} and uniform convergence that
\begin{align*}
    Y^{\vecv}_{t}(\omega)&\leq\inf_{\tau\in\mcT_t}\mcE^{\vecv}_t\Big[\ett_{[\tau=T]}\varphi(\vecv) + \ett_{[\tau<T]}\inf_{b\in U}\{Y^{\vecv\circ(\tau,b),k}_\tau+c(\vecv\circ (\tau,b))\}\Big](\omega)
    \\
    &\leq \inf_{\tau\in\mcT_t}\mcE^{\vecv}_t\Big[\ett_{[\tau=T]}\varphi(\vecv) + \ett_{[\tau<T]}\inf_{b\in U}\{Y^{\vecv\circ(\tau,b)}_\tau+c(\vecv\circ (\tau,b))\}\Big](\omega)+C/k.
\end{align*}
On the other hand, monotonicity implies that
\begin{align*}
    Y^{\vecv,k}_{t}(\omega)&\geq\inf_{\tau\in\mcT_t}\mcE^{\vecv}_t\Big[\ett_{[\tau=T]}\varphi(\vecv) + \ett_{[\tau<T]}\inf_{b\in U}\{Y^{\vecv\circ(\tau,b)}_\tau+c(\vecv\circ (\tau,b))\}\Big](\omega)
\end{align*}
for each $k\geq 0$ and the result follows by letting $k\to\infty$.\qed\\


\section{Value of the game\label{sec:value}}
We show that the game has a value by proving that $Z^{\vecv}_t(\omega)=Y^{\vecv}_t(\omega)$ for all $(t,\omega,\vecv)\in\Lambda\times \mcD$. As a bonus, this immediately gives that the lower value function $Z$ satisfies the DPP in \eqref{ekv:dynp-Y}.

The approach is once again to arrive at the result through a sequence of lemmata. Similarly to the above, we define the truncated lower value function as
\begin{align}
Z^{\vecv,k}_ t(\omega):=\sup_{\Prob^S\in\mcP_t^S}\inf_{u\in\mcU^k_t}\E^{\Prob^S(\vecv\circ u)}_t\Big[\varphi(\vecv\circ u) +\sum_{j=1}^{N} c(\vecv\circ [u]_{j})\Big](\omega)
\end{align}
for all $(t,\omega,\vecv)\in\Lambda\times\mcD$ and have the following:

\begin{lem}\label{lem:Z-unif-conv}
There is a $C>0$ such that
\begin{align*}
Z^{\vecv}_t(\omega)\geq Z^{\vecv,k}_t(\omega)-C/k
\end{align*}
for all $(t,\omega,\vecv)\in \Lambda\times\mcD$ and each $k\geq 0$.
\end{lem}

\noindent \emph{Proof.} We have
\begin{align*}
Z^{\vecv,k}_ t-Z^{\vecv}_ t&=  \sup_{\Prob^S\in\mcP_t^S}\inf_{u\in\mcU^k_t}\E^{\Prob^S(\vecv\circ u)}_t\Big[\varphi(\vecv\circ u) +\sum_{j=1}^{N} c(\vecv\circ [u]_{j})\Big]
\\
&\quad- \sup_{\Prob^S\in\mcP_t^S}\inf_{u\in\mcU_t}\E^{\Prob^S(\vecv\circ u)}_t\Big[\varphi(\vecv\circ u) +\sum_{j=1}^N c(\vecv\circ [u]_{j})\Big]
  \\
  &\leq \sup_{\Prob^S\in\mcP_t^S}\Big\{\inf_{u\in\mcU^k_t}\E^{\Prob^S(\vecv\circ u)}_t\Big[\varphi(\vecv\circ u) +\sum_{j=1}^{N} c(\vecv\circ [u]_{j})\Big]
  \\
  &\quad- \inf_{u\in\mcU_t}\E^{\Prob^S(\vecv\circ u)}_t\Big[\varphi(\vecv\circ u) +\sum_{j=1}^N c(\vecv\circ [u]_{j})\Big]\Big\}.
\end{align*}
Given $\Prob^S\in\mcP^S(t,\omega)$ and $\eps>0$, let $u^\eps=(\tau^\eps_j,\beta^\eps_j)_{j\geq 1}\in\mcU_t$ be such that
\begin{align*}
\inf_{u\in\mcU_t}\E^{\Prob^S(\vecv\circ u)}_t\Big[\varphi(\vecv\circ u) +\sum_{j=1}^{N^\eps} c(\vecv\circ [u]_{j})\Big]\geq \E^{\Prob^S(\vecv\circ u^\eps)}_t\Big[\varphi(\vecv\circ u^\eps) +\sum_{j=1}^{N^\eps} c(\vecv\circ [u^\eps]_{j})\Big]-\eps,
\end{align*}
where $N^\eps:=\max\{j\geq 0:\tau_j<T\}$. Then, since $\varphi$ is bounded and $c\geq\delta>0$,
\begin{align*}
\E^{\Prob^S(\vecv\circ u^\eps)}_t\big[N^\eps\big]\leq \frac{1}{\delta}\E^{\Prob^S(\vecv\circ u^\eps)}_t\Big[\sum_{j=1}^{N^\eps} c(\vecv\circ [u^\eps]_{j})\Big]\leq C.
\end{align*}
Moreover,
\begin{align*}
&\inf_{u\in\mcU^k_t}\E^{\Prob^S(\vecv\circ u)}_t\Big[\varphi(\vecv\circ u) +\sum_{j=1}^{N} c(\vecv\circ [u]_{j})\Big]- \inf_{u\in\mcU_t}\E^{\Prob^S(\vecv\circ u)}_t\Big[\varphi(\vecv\circ u) +\sum_{j=1}^N c(\vecv\circ [u]_{j})\Big]
\\
&\leq \E^{\Prob^S(\vecv\circ [u^\eps]_k)}_t\Big[\varphi(\vecv\circ [u^\eps]_k) +\sum_{j=1}^{N^\eps\wedge k} c(\vecv\circ [u^\eps]_{j})\Big]- \E^{\Prob^S(\vecv\circ u^\eps)}_t\Big[\varphi(\vecv\circ u^\eps) +\sum_{j=1}^{N^\eps}  c(\vecv\circ [u^\eps]_{j})\Big]+\eps.
\end{align*}
Now, since $\Prob^S$ is a non-anticipative map,
\begin{align*}
&\E^{\Prob^S(\vecv\circ [u^\eps]_k)}_t\Big[\varphi(\vecv\circ [u^\eps]_k) +\sum_{j=1}^{N^\eps\wedge k} c(\vecv\circ [u^\eps]_{j})\Big]
\\
&=\E^{\Prob^S(\vecv\circ [u^\eps]_k)}_t\Big[\ett_{[\tau^\eps_{k+1}= T]}\varphi(\vecv\circ u^\eps) +\sum_{j=1}^{N^\eps\wedge k}  c(\vecv\circ [u^\eps]_{j})\Big]+\E^{\Prob^S(\vecv\circ[u^\eps]_k)}_t\Big[\ett_{[\tau^\eps_{k+1}<T]}\varphi(\vecv\circ [u^\eps]_k) \Big]
\\
&\leq\E^{\Prob^S(\vecv\circ u^\eps)}_t\Big[\ett_{[\tau^\eps_{k+1}= T]}\varphi(\vecv\circ u^\eps) +\sum_{j=1}^{N^\eps\wedge k}  c(\vecv\circ [u^\eps]_{j})\Big]+C_0\E^{\Prob^S(\vecv\circ u^\eps)}_t\big[\ett_{[\tau^\eps_{k+1}<T]} \big].
\end{align*}
Put together, we find that
\begin{align*}
&\inf_{u\in\mcU^k_t}\E^{\Prob^S(\vecv\circ u)}_t\Big[\varphi(\vecv\circ u) +\sum_{j=1}^{N} c(\vecv\circ [u]_{j})\Big]- \inf_{u\in\mcU_t}\E^{\Prob^S(\vecv\circ u)}_t\Big[\varphi(\vecv\circ u) +\sum_{j=1}^N c(\vecv\circ [u]_{j})\Big]
\\
&\leq \E^{\Prob^S(\vecv\circ u^\eps)}_t\Big[\ett_{[\tau^\eps_{k+1}= T]}\varphi(\vecv\circ u^\eps)-\varphi(\vecv\circ u^\eps) \Big]+C_0\E^{\Prob^S(\vecv\circ u^\eps)}_t\big[\ett_{[\tau^\eps_{k+1}<T]} \big]+\eps
\\
&\leq 2C_0\E^{\Prob^S(\vecv\circ u^\eps)}_t\big[\ett_{[\tau^\eps_{k+1}<T]} \big]+\eps
\\
&\leq C/(k+1)+\eps.
\end{align*}
As $\eps>0$ was arbitrary and $C$ does not depend on either $\Prob^S$ or $\eps$, the result follows by first letting $\eps\to 0$ and then taking the supremum over all $\Prob^S\in\mcP^S(t,\omega)$.\qed\\

\begin{lem}\label{lem:Zk-prop}
For each $k,\kappa\geq 0$ and $t\in [0,T]$, the map $(\omega,\vecv)\mapsto Z^{\vecv,k}_{t}(\omega):\Omega\times D^\kappa\to\R$ is uniformly continuous with a modulus of continuity that does not depend on $t$.
\end{lem}

\noindent \emph{Proof.} We once again have
\begin{align*}
Z^{\vecv',k}_{t}(\omega')-Z^{\vecv,k}_t(\omega)&\leq\sup_{u\in\mcU^k_{t}}\Big\{\mcE^{\vecv'\circ u}\big[\varphi(\vecv'\circ u) + \sum_{j=1}^N c(\vecv'\circ [u]_j)\big](\omega')\\
&\quad-\mcE^{\vecv\circ u_\omega}\big[\varphi(\vecv\circ u_{\omega}) + \sum_{j=1}^N c(\vecv\circ [u_{\omega}]_j)\big](\omega)\Big\}
\end{align*}
with $u_\omega$ as in Assumption \ref{ass:mcE-consist}. This gives that
\begin{align*}
Z^{\vecv',k}_{t}(\omega')-Z^{\vecv,k}_t(\omega)\leq \rho_{\mcE,\kappa+k}(\|\omega'-\omega\|_t+|\vecv'-\vecv|)
\end{align*}
by again using the same assumption.\qed\\

For $\eps>0$, we introduce the discretization of $Z^{\cdot,k}$ as
\begin{align}
\bar Z^{\vecv,k,\eps}_ t(\omega):=\sup_{\Prob^S\in\mcP_t^S}\inf_{u\in\bar\mcU^{k,\eps}_t}\E^{\Prob^S(\vecv\circ u)}_t\Big[\varphi(\vecv\circ u) +\sum_{j=1}^{N} c(\vecv\circ [u]_{j})\Big](\omega),
\end{align}
where $\bar\mcU^{k,\eps}_t:=\{u\in\mcU^k: (\tau_j(\omega),\beta_j(\omega))\in \bbT^\eps\times\bar U^\eps\text{ for all }\omega\in\Omega\text{ and }j=1,\ldots,k\}$.

\begin{lem}\label{lem:Z-eps-cov}
For each $k,\kappa\geq 0$, there is a modulus of continuity $\rho_{\bar Z,\kappa,k}$, such that
\begin{align}
\bar Z^{\vecv,k,\eps}_t(\omega)\leq Z^{\vecv,k}_ t(\omega)+\rho_{\bar Z,\kappa,k}(\eps),
\end{align}
for all $(t,\omega,\vecv)\in\Lambda\times D^\kappa$.
\end{lem}

\noindent \emph{Proof.} For $u\in\mcU^k$, we let
\begin{align*}
  \Xi^\eps(u):=\Big(\sum_{i=1}^{n^\eps_t}t^\eps_i\ett_{[t^\eps_{i-1},\bar t^\eps_i)}(\tau_j)+T\ett_{[\tau_j=T]}, \sum_{l=1}^{n^\eps_U}b^\eps_l\ett_{U^\eps_l}(\beta_j)\Big)_{j=1}^k
\end{align*}
and note that $\Xi^\eps(u)\in\bar\mcU^{k,\eps}$. We thus have that
\begin{align*}
\bar Z^{\vecv,k,\eps}_ t(\omega)-Z^{\vecv,k}_ t(\omega)&\leq \sup_{\Prob^S\in\mcP_t^S}\sup_{u\in\bar\mcU^{k}_t}\E^{\Prob^S(\vecv\circ u)}_t\Big[|\varphi(\vecv\circ u)-\varphi(\vecv\circ \Xi^\eps(u))
\\
&\quad+\sum_{j=1}^{N} |c(\vecv\circ [u]_{j})-c(\vecv\circ [\Xi^\eps(u)]_{j})|\Big](\omega)
\\ &\leq 2k\sup_{u\in\bar\mcU^{k}_t}\mcE^{\vecv\circ u}_t\Big[(\rho_{\varphi,\kappa+k}+\rho_{\varphi,\kappa+k})(\textbf d_{\kappa+k}[(T,\cdot,\vecv\circ u),(T,\cdot,\vecv\circ \Xi^n(u))])\Big](\omega)
\\ &\leq \rho(\eps)
\end{align*}
for some modulus of continuity $\rho$.\qed\\

\begin{lem}\label{lem:dynp-bar-Z}
For each $\eps > 0$, the family $(\bar Z^{\cdot,k,\eps})_{k\geq 0}$ satisfies the recursion
\begin{align}\label{ekv:dynp-bar-Z}
\begin{cases}
  \bar Z^{\vecv,k,\eps}_{T}(\omega)&=\varphi(\vecv)\\
  \bar Z^{\vecv,k,\eps}_{t^\eps_i}(\omega)&=\min_{b\in \bar U^\eps}\big\{\bar Z^{\vecv\circ(t^\eps_i,b),k-1,\eps}_{t^\eps_i}+c(\vecv\circ(t^\eps_i,b))\big\}\wedge \mcE^\vecv_{t^\eps_i}\big[\bar Z^{\vecv,k,\eps}_{t^\eps_{i+1}}\big](\omega),\quad i=1,\ldots,n^\eps_t-1
\end{cases}
\end{align}
for all $k\geq 1$. Moreover,
\begin{align}\label{ekv:Z-k-eps-is-Y-k-eps}
\bar Z^{\vecv,k,\eps}_{t^\eps_i}(\omega)=\inf_{u\in\bar\mcU^{k,\eps}_{t^\eps_i}}\mcE^{\vecv\circ u}_{t^\eps_i}\Big[\varphi(\vecv\circ u) +\sum_{j=1}^{N} c(\vecv\circ [u]_{j})\Big](\omega)=:\bar Y^{\vecv,k,\eps}_{t^\eps_i}(\omega)
\end{align}
for all $k\geq 0$, $i=1,\ldots,n^\eps_t$ and $(\omega,\vecv)\in\Omega\times\mcD$.
\end{lem}

\noindent \emph{Proof.} The result will follow by a double induction and we note that the terminal condition holds trivially since \eqref{ekv:no-imp@end} implies that
\begin{align}\label{ekv:equal@end}
  \bar Z^{\vecv,k,\eps}_{T}(\omega)=\bar Y^{\vecv,k,\eps}_{T}(\omega)=\varphi(\vecv)
\end{align}
for all $(\omega,\vecv)\in\Omega\times \mcD$ and $k\geq 0$. We let $(R^{\cdot,k})_{k\geq 0}$ be such that $R^{\cdot,0}\equiv \bar Z^{\cdot,0,\eps}$ and $(R^{\cdot,k})_{k\geq 1}$ is the unique solution to the recursion \eqref{ekv:dynp-bar-Z} and assume that for some $k'\geq 1$ and $i'\in \{1,\ldots,n_t^\eps-1\}$ the statement holds for all $i\in \{1,\ldots,n_t^\eps-1\}$ when $k=1,\ldots,k'-1$ and for $i=i'+1,\ldots,n^\eps_t$ when $k=k'$.\\

\textbf{Step 1.} We first show that $\bar Y^{\vecv,k,\eps}_{t^\eps_{i'}}$ satisfies the recursion step in \eqref{ekv:dynp-bar-Z}. By definition we have,
\begin{align*}
\bar Y^{\vecv,k,\eps}_{t^\eps_{i'}}&=  \inf_{u\in\bar\mcU^{k,\eps}_{t^\eps_{i'}}}(\ett_{[\tau_1=t^\eps_{i'}]}+\ett_{[\tau_1>t^\eps_{i'}]})\mcE^{\vecv\circ u}_{t^\eps_{i'}}\Big[\varphi(\vecv\circ u) +\sum_{j=1}^{N} c(\vecv\circ [u]_{j})\Big]
\\
&= \inf_{u\in\bar\mcU^{k-1,\eps}_{t^\eps_{i'}}}\min_{b\in \bar U^\eps}\Big\{\mcE^{\vecv\circ (t^\eps_{i'},b)\circ u}_{t^\eps_{i'}}\Big[\varphi(\vecv\circ(t^\eps_i,b)\circ u) +\sum_{j=1}^{N} c(\vecv\circ (t^\eps_{i'},b)\circ [u]_{j})\Big]+c(\vecv\circ(t^\eps_i,b))\Big\}
\\
&\quad\wedge\inf_{u\in\bar\mcU^{k,\eps}_{t^\eps_{i'+1}}}\mcE^{\vecv\circ u}_{t^\eps_{i'}}\Big[\varphi(\vecv\circ u) +\sum_{j=1}^{N} c(\vecv\circ [u]_{j})\Big]
\\
&= \min_{b\in \bar U^\eps}\big\{\bar Y^{\vecv\circ(t^\eps_i,b),k-1,\eps}_{t^\eps_i}+c(\vecv\circ(t^\eps_i,b))\big\}\wedge \inf_{u\in\bar\mcU^{k,\eps}_{t^\eps_{i'+1}}}\mcE^{\vecv\circ u}_{t^\eps_{i'}}\Big[\varphi(\vecv\circ u) +\sum_{j=1}^{N} c(\vecv\circ [u]_{j})\Big].
\end{align*}
By the tower property we have
\begin{align*}
\inf_{u\in\bar\mcU^{k,\eps}_{t^\eps_{i'+1}}}\mcE^{\vecv}_{t^\eps_{i'}}\Big[\varphi(\vecv\circ u) +\sum_{j=1}^{N} c(\vecv\circ [u]_{j})\Big]&=\inf_{u\in\bar\mcU^{k,\eps}_{t^\eps_{i'+1}}}\mcE^{\vecv}_{t^\eps_{i'}}\Big[\mcE^{\vecv\circ u}_{t^\eps_{i'+1}}\Big[\varphi(\vecv\circ u) +\sum_{j=1}^{N} c(\vecv\circ [u]_{j})\Big]\Big]
\\
&\geq \mcE^{\vecv}_{t^\eps_{i'}}\Big[\inf_{u\in\bar\mcU^{k,\eps}_{t^\eps_{i'+1}}}\mcE^{\vecv\circ u}_{t^\eps_{i'+1}}\Big[\varphi(\vecv\circ u) +\sum_{j=1}^{N} c(\vecv\circ [u]_{j})\Big]\Big]
\\
&=\mcE^{\vecv}_{t^\eps_{i'}}\big[\bar Y^{\vecv,k,\eps}_{t^\eps_{i'+1}}\big],
\end{align*}
establishing that $\bar Y^{\vecv,k,\eps}_{t^\eps_{i'}}\geq R^{\vecv,k}_{t^\eps_{i'}}$. To arrive at the opposite inequality we introduce the impulse control $u^*:=(\tau_j^*,\beta^*_j)_{j=1}^{N^*}$, where
\begin{align*}
  \tau_j^*:=\min\big\{s\in\bbT^\eps_{\tau^*_{j-1}}: R^{\vecv\circ [u^*]_{j-1},k+1-j}_{s}=\min_{b\in \bar U^\eps}\{R^{\vecv\circ [u^*]_{j-1}\circ(s,b),k-1}_{s}+c(\vecv\circ [u^*]_{j-1}\circ(s,b))\}\big\}
\end{align*}
for $j=1\ldots,k$, with $\tau^*_0:=t_{i'}$ and $\bbT^\eps_{\tau}(\omega'):=\{s\in \bbT^\eps:s\geq \tau(\omega')\}$, $N^*:=\min\{j\geq 0: \tau_j^*<T\}\wedge k$ and $\beta^*_j$ is measurable selection of
\begin{align*}
  \beta_j^*\in\argmin_{b\in \bar U^\eps}\big\{R^{\vecv\circ [u^*]_{j-1}\circ(\tau_j^*,b),k-1}_{\tau_j^*}+c(\vecv\circ [u^*]_{j-1}\circ(\tau_j^*,b))\big\}.
\end{align*}
Then, $u^*\in\bar\mcU^{k,\eps}_{t_{i'}}$ and it easily follows by repeated use of \eqref{ekv:dynp-bar-Z} and the tower property that
\begin{align*}
  R^{\vecv,k}_{t^\eps_{i'}}=\mcE^{\vecv\circ u^*}_{t^\eps_{i'}}\Big[\varphi(\vecv\circ u^*) +\sum_{j=1}^{N} c(\vecv\circ [u^*]_{j})\Big]\geq \bar Y^{\vecv,k,\eps}_{t^\eps_{i'}}
\end{align*}
and we conclude that $\bar Y^{\vecv,k,\eps}_{t^\eps_{i'}}$ satisfies the recursion step in \eqref{ekv:dynp-bar-Z}.\\

\textbf{Step 2.} For the sake of completeness, we also prove rigorously that $\bar Z^{\vecv,k,\eps}_{t^\eps_i}\leq \bar Y^{\vecv,k,\eps}_{t^\eps_i}$. We have
\begin{align*}
\bar Z^{\vecv,k,\eps}_{t^\eps_{i'}}&=  \sup_{\Prob^S\in\mcP_{t^\eps_{i'}}^S}\inf_{u\in\bar\mcU^{k,\eps}_{t^\eps_{i'}}}(\ett_{[\tau_1=t^\eps_{i'}]}+\ett_{[\tau_1>t^\eps_{i'}]})\E^{\Prob^S(\vecv\circ u)}_{t^\eps_{i'}}\Big[\varphi(\vecv\circ u) +\sum_{j=1}^{N} c(\vecv\circ [u]_{j})\Big]
\\
&= \sup_{\Prob^S\in\mcP_{t^\eps_{i'}}^S}\Big\{\inf_{u\in\bar\mcU^{k-1,\eps}_{t^\eps_{i'}}}\min_{b\in \bar U^\eps}\Big\{\E^{\Prob^S(\vecv\circ (t^\eps_{i'},b)\circ u)}_{t^\eps_{i'}}\Big[\varphi(\vecv\circ(t^\eps_i,b)\circ u) +\sum_{j=1}^{N} c(\vecv\circ (t^\eps_{i'},b)\circ [u]_{j})\Big]+c(\vecv\circ(t^\eps_i,b))\Big\}
\\
&\quad\wedge\inf_{u\in\bar\mcU^{k,\eps}_{t^\eps_{i'+1}}}\E^{\Prob^S(\vecv\circ u)}_{t^\eps_{i'}}\Big[\varphi(\vecv\circ u) +\sum_{j=1}^{N} c(\vecv\circ [u]_{j})\Big]\Big\}
\\
&= \min_{b\in \bar U^\eps}\big\{\bar Z^{\vecv\circ(t^\eps_i,b),k-1,\eps}_{t^\eps_i}+c(\vecv\circ(t^\eps_i,b))\big\}\wedge \sup_{\Prob^S\in\mcP_{t^\eps_{i'}}^S}\inf_{u\in\bar\mcU^{k,\eps}_{t^\eps_{i'+1}}}\E^{\Prob^S(\vecv\circ u)}_{t^\eps_{i'}}\Big[\varphi(\vecv\circ u) +\sum_{j=1}^{N} c(\vecv\circ [u]_{j})\Big].
\end{align*}
By the tower property for $\E^\Prob$ and the induction hypothesis we have
\begin{align*}
&\sup_{\Prob^S\in\mcP_{t^\eps_{i'}}^S}\inf_{u\in\bar\mcU^{k,\eps}_{t^\eps_{i'+1}}}\E^{\Prob^S(\vecv\circ u)}_{t^\eps_{i'}}\Big[\varphi(\vecv\circ u) +\sum_{j=1}^{N} c(\vecv\circ [u]_{j})\Big]
\\
&=\sup_{\Prob^S\in\mcP_{t^\eps_{i'}}^S}\E^{\Prob^S(\vecv)}_{t^\eps_{i'}}\Big[\inf_{u\in\bar\mcU^{k,\eps}_{t^\eps_{i'+1}}} \E^{\Prob^S(\vecv\circ u)}_{t^\eps_{i'+1}}\Big[\varphi(\vecv\circ u) +\sum_{j=1}^{N} c(\vecv\circ [u]_{j})\Big]\Big]
\\
&\leq \mcE^{\vecv}_{t^\eps_{i'}}\Big[\inf_{u\in\bar\mcU^{k,\eps}_{t^\eps_{i'+1}}} \mcE^{\vecv\circ u}_{t^\eps_{i'+1}}\Big[\varphi(\vecv\circ u) +\sum_{j=1}^{N} c(\vecv\circ [u]_{j})\Big]\Big]=\mcE^{\vecv}_{t^\eps_{i'}}\big[\bar Z^{\vecv,k,\eps}_{t^\eps_{i'+1}}\big]
\end{align*}
and we conclude that
\begin{align*}
\bar Z^{\vecv,k,\eps}_{t^\eps_{i'}}\leq \hat Z^{\vecv,k,\eps}_{t^\eps_{i'}}:=\min_{b\in \bar U^\eps}\{\bar Z^{\vecv\circ(t^\eps_{i'},b),k-1,\eps}_{t^\eps_{i'}}+c(\vecv\circ(t^\eps_i,b))\}\wedge \mcE^\vecv_{t^\eps_{i'}}\big[\bar Z^{\vecv,k,\eps}_{t^\eps_{{i'}+1}}\big].
\end{align*}
In particular, our induction assumption then implies that $\bar Z^{\vecv,k,\eps}_{t^\eps_{i'}}(\omega)\leq \bar Y^{\vecv,k,\eps}_{t^\eps_i}(\omega)$.\\

\textbf{Step 3.} Finally, we show that $\bar Z^{\vecv,k,\eps}_{t^\eps_{i'}}\geq \bar Y^{\vecv,k,\eps}_{t^\eps_{i'}}$. We do this by showing that for any $\eps'>0$, there is a strategy $\Prob^{S,\eps'}\in\mcP^S({t^\eps_{i'}},\omega)$ such that
\begin{align*}
  \bar Y^{\vecv,k,\eps}_{t^\eps_{i'}}(\omega)\leq \E^{\Prob^{S,\eps'}(\vecv\circ u)}_t\Big[\varphi(\vecv\circ u) +\sum_{j=1}^{N} c(\vecv\circ [u]_{j})\Big](\omega)+\eps'
\end{align*}
for each $u\in \bar\mcU^{k,\eps}_{t^\eps_{i'}}$ and $\omega\in\Omega$. First note that for any $u\in \bar\mcU^{k,\eps}_{t^\eps_{i'}}$ and $i\in\{i'+1,\ldots,n^\eps_t\}$, the map
\begin{align*}
\omega'\mapsto \bar Y^{\vecv\circ u_{t^\eps_i}(\omega\otimes_{t^\eps_{i'}}\omega'),k-N(t^\eps_i)(\omega\otimes_{t^\eps_{i'}}\omega'),\eps}_{t^\eps_i}(\omega\otimes_{t^\eps_{i'}}\omega')
\end{align*}
is Borel-measurable and thus upper semi-analytic.  Given $u\in \bar\mcU^{k,\eps}_{t^\eps_{i'}}$ we can then, by repeatedly arguing as in the proof of Lemma 4.11 in \cite{NutzZhang15}, find a sequence of measures $(\Prob^u_{i})_{i=i'}^{n^\eps_t-1}$ with $\Prob^u_{i}\in\mcP(t^\eps_{i'},\omega,\vecv\circ u_{t^\eps_{i'}})$ such that $\Prob^u_{i}=\Prob^u_{i-1}$ on $\mcF_{t^\eps_{i-1}-t^\eps_{i'}}$ for $i=i'+1,\ldots,n^\eps_{t}-1$ and
\begin{align*}
  \E^{\Prob^u_{i}}\Big[\bar Y^{\vecv\circ u_{t^\eps_{i}}(\omega\otimes_{t^\eps_{i'}}\cdot),k-N(t^\eps_{i})(\omega\otimes_{t^\eps_{i'}}\cdot),\eps}_{t^\eps_{i+1}}(\omega\otimes_{t^\eps_{i'}}\cdot)\Big| \mcF_{t^\eps_i-t^\eps_{i'}}\Big]\geq \mcE^{\vecv\circ u}_{t^\eps_{i}}\big[\bar Y^{\vecv\circ u_{t^\eps_{i}},k-N(t^\eps_i),\eps}_{t^\eps_{i+1}}](\omega\otimes_{t^\eps_{i'}}\cdot)-\eps'/n^\eps_t.
\end{align*}
This implies the existence of a strategy $\Prob^{S,\eps'}\in\mcP^S({t^\eps_{i'}},\omega)$ by setting $\Prob^{S,\eps'}(\vecv\circ u):=\Prob^u_{n^{\eps}_t}$ such that for any $u\in\bar\mcU^k_{t^\eps_{i'}}$, we have
\begin{align*}
  \bar Y^{\vecv,k,\eps}_{t^\eps_{i'}}&\leq\bar Y^{\vecv\circ u_{t^\eps_{i'}},k-N(t^{\eps}_{i'}),\eps}_{t^\eps_{i'}}+\sum_{j=1}^{N(t^{\eps}_{i'})}c(\vecv\circ[u]_j)
  \\
  &\leq \E^{\Prob^{S,\eps'}(\vecv\circ u)}_{t^\eps_{i'}}\Big[\bar Y^{\vecv\circ u_{t^\eps_{i'}},k-N(t^{\eps}_{i'}),\eps}_{t^\eps_{i'+1}}\Big]+\sum_{j=1}^{N(t^{\eps}_{i'})}c(\vecv\circ[u]_j)+\eps'/n^{\eps}_t
  \\
  &\leq \E^{\Prob^{S,\eps'}(\vecv\circ u)}_{t^\eps_{i'}}\Big[\bar Y^{\vecv\circ u_{t^\eps_{i'+1}},k-N(t^{\eps}_{i'+1}),\eps}_{t^\eps_{i'+1}}+\sum_{j=N(t^{\eps}_{i'})+1}^{N(t^{\eps}_{i'+1})}c(\vecv\circ[u]_j)\Big]+\sum_{j=1}^{N(t^{\eps}_{i'})}c(\vecv\circ[u]_j)+\eps'/n^{\eps}_t
  \\
  &=\E^{\Prob^{S,\eps'}(\vecv\circ u)}_{t^\eps_{i'}}\Big[\bar Y^{\vecv\circ u_{t^\eps_{i'+1}},k-N(t^{\eps}_{i'+1}),\eps}_{t^\eps_{i'+1}}+\sum_{j=1}^{N(t^{\eps}_{i'+1})}c(\vecv\circ[u]_j)\Big]+\eps'/n^{\eps}_t.
\end{align*}
By repeating this process $n^\eps_t-i'$ times followed by using \eqref{ekv:equal@end} and the tower property for $\E^{\Prob^{S,\eps'}(\vecv\circ u)}$ we find that
\begin{align*}
  \bar Y^{\vecv,k,\eps}_{t^\eps_i}(\omega)\leq \E^{\Prob^{S,\eps'}(\vecv\circ u)}_{t^\eps_i}\Big[\varphi(\vecv\circ u) +\sum_{j=1}^{N} c(\vecv\circ [u]_{j})\Big](\omega)+\eps'.
\end{align*}
and since  $u\in\bar\mcU^k_{t^\eps_{i'}}$ was arbitrary, we conclude that $\bar Y^{\vecv,k,\eps}_{t^\eps_{i'}}(\omega)\leq \bar Z^{\vecv,k,\eps}_{t^\eps_{i'}}(\omega)+\eps'$ from which the assertion follows as $\eps'>0$ was arbitrary. This concludes the induction step.\qed\\

\begin{cor}
We can extend the definition of $\bar Y$ to $\Lambda\times\mcD$ by letting
\begin{align*}
  \bar Y^{\vecv,k,\eps}_{t}(\omega):=\inf_{u\in\bar\mcU^{k,\eps}_{t}}\mcE^{\vecv\circ u}_{t}\Big[\varphi(\vecv\circ u) +\sum_{j=1}^{N} c(\vecv\circ [u]_{j})\Big](\omega)
\end{align*}
for all $(t,\omega,\vecv)\in\Lambda\times\mcD$ and have that $\bar Y^{\cdot,k,\eps}\equiv \bar Z^{\cdot,k,\eps}$ for all $k\geq 0$.
\end{cor}

\noindent \emph{Proof.} Repeating the above proof we find that $\bar Y^{\cdot,k,\eps}$ and $\bar Z^{\cdot,k,\eps}$ both satisfy the recursion
\begin{align*}
  \bar Z^{\vecv,k,\eps}_{t}(\omega)&=\mcE^\vecv_{t}\big[\bar Z^{\vecv,k,\eps}_{t^\eps_{i+1}}\big](\omega)
\end{align*}
for all $t\in (t^\eps_i,t^\eps_{i+1})$.\qed\\

\begin{thm}\label{thm:game-value}
The game has a value, \ie $Y\equiv Z$.
\end{thm}

\noindent \emph{Proof.} Applying an argument identical to the one in the proof of Lemma~\ref{lem:Z-eps-cov} gives that $\|\bar Y^{\vecv,k,\eps}-Y^{\vecv,k}\|_T\to 0$ as $\eps\to 0$. By uniqueness of limits we find that $Y^{\cdot,k}=Z^{\cdot,k}$ and taking the limit as $k$ tends to infinity, the result follows by Lemma~\ref{lem:unif-conv} and Lemma~\ref{lem:Z-unif-conv}.\qed\\

\begin{cor}\label{thm:dynP-Z}
The map $Z$ is bounded, uniformly continuous and satisfies the recursion
\begin{align}\label{ekv:dynp-Z}
    Z^{\vecv}_t(\omega)=\inf_{\tau\in\mcT_t}\mcE^{\vecv}_t\Big[\ett_{[\tau=T]}\varphi(\vecv) + \ett_{[\tau<T]}\inf_{b\in U}\{Z^{\vecv\circ(\tau,b)}_\tau+c^{}(\vecv\circ(\tau,b))\}\Big](\omega),
  \end{align}
for all $(t,\omega,\vecv)\in\Lambda\times\mcD$.
\end{cor}

\noindent \emph{Proof.} Since $Z\equiv Y$ the statement in Theorem~\ref{thm:dynP} applies to $Z$ as well.\qed\\

\section{A verification theorem\label{sec:verif-thm}}
In the previous two sections we have shown that the dynamic programming principle is a necessary condition for a map to be the value function of \eqref{ekv:minmax-form}: if $\mcY$ is an upper or lower value function then it satisfies \eqref{ekv:dynp-Y}. In this section we turn to sufficiency and its implications, starting with the following uniqueness result:
\begin{prop}\label{prop:uniqueness}
Suppose that there is a progressively measurable map $\mcY:\Lambda\times \mcD\to\R$ that is uniformly continuous and bounded, such that
\begin{align}
    \mcY^{\vecv}_t(\omega)=\inf_{\tau\in\mcT_t}\mcE^{\vecv}_t\Big[\ett_{[\tau=T]}\varphi(\vecv) + \ett_{[\tau<T]}\inf_{b\in U}\{\mcY^{\vecv\circ(\tau,b)}_\tau+c(\vecv\circ(\tau,b))\}\Big](\omega)
\end{align}
for all $(t,\omega,\vecv)\in\Lambda\times\mcD$. Then $\mcY\equiv Y$.
\end{prop}

\noindent \emph{Proof.} For any $\hat u=(\hat\tau_j,\hat\beta_j)_{j=1}^k\in\mcU^k_t$ we have
\begin{align*}
 \mcY^\vecv_t\leq\mcE^{\vecv}_t\Big[\ett_{[\hat\tau_1=T]}\varphi(\vecv) + \ett_{[\hat\tau_1<T]}\{\mcY^{\vecv\circ(\hat\tau_1,\hat\beta_1)}_{\hat\tau_1}+c(\vecv\circ(\hat\tau_1,\hat\beta_1))\}\Big].
\end{align*}
But similarly,
\begin{align*}
  \mcY^{\vecv\circ(\hat\tau_1,\hat\beta_1)}_{\hat\tau_1}&\leq \mcE^{\vecv\circ(\hat\tau_1,\hat\beta_1)}_{\hat\tau_1}\Big[\ett_{[\hat\tau_2=T]}\varphi(\vecv\circ(\hat\tau_1,\hat\beta_1)) + \ett_{[\hat\tau_2<T]}\{\mcY^{\vecv\circ[\hat u]_2}_{\tau^\eps_2}+c(\vecv\circ[\hat u]_2)\}\Big].
\end{align*}
Repeating this argument $k$ times we find, since $\hat u\in\mcU^k_t$ was arbitrary, that
\begin{align*}
  \mcY^\vecv_t(\omega)&\leq\inf_{u\in\mcU^k_t}\mcE^{\vecv\circ u}_t\big[\ett_{[\tau_{k}=T]}\varphi(\vecv\circ u) + \sum_{j=1}^{N} c(\vecv\circ [u]_j)+\ett_{[\tau_k<T]}\{\mcY^{\vecv\circ u}_{\tau_k}+c(\vecv\circ u)\}\big]
  \\
  &\leq Y^{\vecv,k-1}_t(\omega).
\end{align*}
Letting $k\to\infty$, it follows by Lemma~\ref{lem:unif-conv} that $\mcY^\vecv_t(\omega)\leq Y^\vecv_t(\omega)$ for all $(t,\omega,\vecv)\in \Lambda\times\mcD$.

On the other hand, from Theorem~\ref{thm:NutzZhang} and Lemma~\ref{lem:eps-attaind} there is for each $\eps>0$ a pair $(\tau^\eps_1,\beta^\eps_1)\in\mcU^1_t$ such that
\begin{align*}
  \mcY^\vecv_t\geq\mcE^{\vecv}_t\Big[\ett_{[\tau^\eps_1=T]}\varphi(\vecv) + \ett_{[\tau^\eps_1<T]}\{\mcY^{\vecv\circ(\tau^\eps_1,\beta^\eps_1)}_{\tau^\eps_1}+c(\vecv\circ(\tau^\eps_1,\beta^\eps_1))\}\Big]-\eps/2.
\end{align*}
Moreover, Step 1 in the proof of Lemma~\ref{lem:dynp-Yk} implies the existence of a $(\tau^\eps_2,\beta^\eps_2)\in\mcU^1_{\tau^\eps_1}$ such that
\begin{align*}
  \mcY^{\vecv\circ(\tau^\eps_1,\beta^\eps_1)}_{\tau^\eps_1}&\geq \mcE^{\vecv\circ(\tau^\eps_1,\beta^\eps_1)}_t\Big[\ett_{[\tau^\eps_2=T]}\varphi(\vecv\circ(\tau^\eps_1,\beta^\eps_1))
  \\
  &\quad+ \ett_{[\tau^\eps_2<T]}\{\mcY^{\vecv\circ(\tau^\eps_1,\beta^\eps_1)\circ(\tau^\eps_2,\beta^\eps_2)}_{\tau^\eps_2}+c(\vecv\circ(\tau^\eps_1,\beta^\eps_1)\circ(\tau^\eps_2,\beta^\eps_2))\}\Big]-\eps/4.
\end{align*}
Repeating this argument indefinitely gives us an infinite sequence $u^\eps=(\tau_j^\eps,\beta_j^\eps)_{j= 1}^\infty\in\mcU_t$ such that
\begin{align*}
  \mcY^{\vecv\circ [u^\eps]_{k-1}}_{\tau^\eps_{k-1}}&\geq \mcE^{\vecv\circ [u^\eps]_{k-1}}_{\tau^\eps_{k-1}}\Big[\ett_{[\tau^\eps_k=T]}\varphi(\vecv\circ [u^\eps]_{k-1}) + \ett_{[\tau^\eps_k<T]}\{\mcY^{\vecv\circ [u^\eps]_{k})}_{\tau^\eps_2}+c(\vecv\circ [u^\eps]_{k})\}\Big]-\eps/2^k,
\end{align*}
for each $k\geq 1$ and we find by applying the tower property that
\begin{align*}
  \mcY^\vecv_t\geq\mcE^{\vecv\circ u^\eps}_t\big[\ett_{[\tau^\eps_{k}=T]}\varphi(\vecv\circ u^\eps) + \sum_{j=1}^{N^\eps\wedge k} c([u^\vecv]_j)+\ett_{[\tau^\eps_k<T]}\mcY^{\vecv\circ [u^\eps]_k}_{\tau^\eps_k}\big]-\eps.
\end{align*}
Taking the limit as $k\to\infty$ on the right hand side thus gives
\begin{align*}
  \mcY^\vecv_t+\eps&\geq\liminf_{k\to\infty}\sup_{\Prob\in\mcP(t,\omega,\vecv\circ u^\eps)}\E^{\Prob}_t\big[\ett_{[\tau^\eps_{k}=T]}\varphi(\vecv\circ u^\eps) + \sum_{j=1}^{N^\eps\wedge k} c(\vecv\circ [u^\eps]_j)+\ett_{[\tau^\eps_k<T]}\mcY^{\vecv\circ[ u^\eps]_k}_{\tau^\eps_k}\big]
  \\
  &\geq \sup_{\Prob\in\mcP(t,\omega,\vecv\circ u^\eps)}\liminf_{k\to\infty}\E^{\Prob}_t\big[\ett_{[\tau^\eps_{k}=T]}\varphi(\vecv\circ u^\eps) + \sum_{j=1}^{N^\eps\wedge k} c(\vecv\circ [u^\eps]_j)+\ett_{[\tau^\eps_k<T]}\mcY^{\vecv\circ [u^\eps]_k}_{\tau^\eps_k}\big].
\end{align*}
Now, since $\mcY$ and $\varphi$ are uniformly bounded and $c\geq\delta>0$, arguing as in the proof of Lemma~\ref{lem:unif-conv} gives that $\mcE^{\vecv\circ u^\eps}\big[\ett_{[\tau^\eps_k<T]}\big]\leq C/k$. In particular, $\Prob\big[\{\tau^\eps_k<T,\,\forall k\geq 0\}\big]= 0$ for any $\Prob\in\mcP(t,\omega,\vecv\circ u^\eps)$ and we can apply Fatou's lemma in the above inequality to conclude that
\begin{align*}
  \mcY^\vecv_t+\eps&\geq\mcE^{\vecv\circ u^\eps}_t\big[\varphi(\vecv\circ u^\eps) + \sum_{j=1}^{N^\eps} c(\vecv\circ [u^\eps]_j)\big]\geq Y^{\vecv} _t.
\end{align*}
Since $\eps>0$ was arbitrary, it follows that $\mcY\equiv Y$.\qed\\

Having proved that the value function is the unique solution to the dynamic programming equation \eqref{ekv:dynp-Y} we show that, under additional measurability and compactness assumptions, $Y$ can be used to extract an optimal control/strategy pair.

\begin{thm}\label{thm:optim-ctrls}
Assume that $u^*=(\tau_j^*,\beta_j^*)_{j=1}^{\infty}$ is such that:
\begin{itemize}
  \item the sequence $(\tau^*_j)_{j=1}^{\infty}$ is given by
\begin{align}
  \tau^*_j:=\inf \Big\{&s \geq \tau^*_{j-1}:\:Y_s^{[u^*]_{j-1}}=\inf_{b\in U} \big\{Y^{[u^*]_{j-1}\circ(s,b)}_s+c([u^*]_{j-1}\circ(s,b))\big\}\Big\}\wedge T,\label{ekv:taujDEF}
\end{align}
using the convention that $\inf\emptyset=\infty$, with $\tau^*_0=0$;
  \item the sequence $(\beta_j^*)_{j=1}^{\infty}$ is such that $\beta^*_j$ is $\mcF_{\tau_j^*}$-measurable and satisfies
\begin{equation}\label{ekv:betajDEF}
  \beta^*_j\in\mathop{\arg\min}_{b\in U}\big\{Y^{[u^*]_{j-1}\circ(\tau^*_j,b)}_{\tau^*_j}+c([u^*]_{j-1}\circ(\tau^*_j,b)\big\}.
\end{equation}
\end{itemize}
Then, $u^*\in\mcU$ is an optimal impulse control for \eqref{ekv:minmax-form} in the sense that
\begin{align}\label{ekv:ustar-is-opt}
  Y^{\emptyset}_0=\mcE^{u^*}\Big[\varphi(u^*) + \sum_{j=1}^{N^*} c([u^*]_j)\Big].
\end{align}
Moreover, if $\mcP(t,\omega,\vecv)$ is weekly compact for each $(t,\omega,\vecv)\in\Lambda\times\mcD$ and $\mcP^{\Prob}_{\tau}(u):=\{\Prob'\in\mcP(u):\Prob'=\Prob\text{ on }\mcF_{\tau}\}$ is weekly compact for all $u\in\mcU$, $\Prob\in\mcP(u)$ and $\tau\in\mcT$, then there is an optimal response $\Prob^{*,S}\in\mcP^S$ for which
\begin{align}\label{ekv:Pstar-is-opt}
  \mcE^{u^*}\big[\varphi(u^*) + \sum_{j=1}^{N^*} c([u^*]_j)\big]=\E^{\Prob^{*,S}(u^*)}\big[\varphi(u^*) + \sum_{j=1}^{N^*} c([u^*]_j)\big]=\inf_{u\in\mcU}\E^{\Prob^{*,S}(u)}\big[\varphi(u) + \sum_{j=1}^{N} c([u]_j)\big].
\end{align}
\end{thm}

\noindent \emph{Proof of \eqref{ekv:ustar-is-opt}.} Repeated use of the dynamic programming principle for $Y$ gives that
\begin{align*}
  Y^{\emptyset}_0&=\mcE^{u^*}\big[\ett_{[\tau_{k}^*=T]}\varphi(u^*) + \sum_{j=1}^{N^*\wedge k} c([u^*]_j)+\ett_{[\tau^*_k<T]}Y^{[u^*]_k}_{\tau^*_k}\big].
\end{align*}
Taking the limit as $k\to\infty$ and repeating the argument in the second part of the proof of Proposition \ref{prop:uniqueness}, \eqref{ekv:ustar-is-opt} follows.\qed\\

To prove the second statement we need two lemmas, the first of which (loosely speaking) shows that for any measure/control pair $(\Prob,u)$ with $u\in\mcU^k$ and $\Prob\in\mcP(u)$, we can extend $\Prob$ optimally from $\tau_k$ until the time that $(Y^u_t)_{\tau_k\leq t\leq T}$ hits the corresponding barrier. This extension then acts as a minimizer up until the first hitting time in \eqref{ekv:dynp-Y}. The main obstacle we face is that $Y^u$ is not necessarily continuous in $\omega$, disqualifying direct use of Lemma 4.5 in \cite{EkrenStopping} as in, for example, Lemma 4.13 of \cite{NutzZhang15}.

\begin{lem}\label{lem:has-P-kernel-1}
For $k\geq 1$ and $u\in\mcU^k$ we let
\begin{align*}
  \tau^\diamond:=\inf \Big\{&s \geq \tau_{k}:\:Y_s^{u}=\inf_{b\in U} \{Y^{u\circ(s,b)}_s+c(u\circ(s,b))\}\Big\}\wedge T
\end{align*}
and assume that for some $\Prob\in \mcP(u)$ the set $\mcP^{\Prob}_{\tau_k}(u)$ is weakly compact. Then, there is a $\Prob^\diamond\in \mcP^{\Prob}_{\tau_k}(u)$ such that
\begin{align}
&\E^{\Prob^\diamond}\Big[\ett_{[\tau^\diamond=T]}\varphi(u) + \ett_{[\tau^\diamond<T]}\inf_{b\in U}\{Y^{u\circ(\tau^\diamond,b)}_{\tau^\diamond}+c(u\circ(\tau^\diamond,b))\}\Big]=\E^{\Prob}\big[Y^{u}_{\tau_k}\big].\label{ekv:P-diamond-attains}
\end{align}
\end{lem}

\noindent \emph{Proof.} To simplify notation we let
\begin{align*}
  S^u_t:=\ett_{[t=T]}\varphi(u) + \ett_{[t<T]}\inf_{b\in U}\{Y^{u\circ(t,b)}_t+c(u\circ(t,b))\}
\end{align*}
and consider the following sequence of stopping times $\eta_l:=\inf\{s\geq \eta_{l-1}: S^{u}_s-Y^{u}_s\leq 1/l\}\wedge T$ for $l\geq 1$ with $\eta_0:=\tau_{k}$. Then, by Proposition 7.50.b) of \cite{BertsekasShreve} and a standard approximation result there is for each $l\geq 0$ a $\mcF_{\tau_k}$-measurable kernel $\nu^l:\Omega\to\PrM(\Omega)$ such that $\nu^l(\cdot)\in \mcP(\tau_k,\cdot,u)$ and
\begin{align}\label{ekv:prob_l-def}
\E^{\nu^l(\omega)}_{\tau_k(\omega)}\big[S^{u}_{\eta_l}\big](\omega)\geq\mcE^{u}_{\tau_k}\big[S^{u}_{\eta_l}\big](\omega)-2^{-l}
\end{align}
for $\Prob$-a.e.~$\omega\in\Omega$. We then define the measure $\Prob_l\in\mcP(u)$ as
\begin{align*}
    \Prob_l[A]:=\int\!\!\!\int (\ett_A)^{\tau_k,\tilde \omega}(\omega')\nu^l(d\omega';\tilde\omega)\Prob[d\tilde\omega]
\end{align*}
for each $A\in\mcF$. Since $\mcP^{\Prob}_{\tau_k}(u)$ is weakly compact we may assume, by possible going to a subsequence, that $\Prob_l\to \Prob^\diamond$ weakly for some $\Prob^\diamond\in \mcP^{\Prob}_{\tau_k}(u)$. We need to show that $\Prob^\diamond$ satisfies \eqref{ekv:P-diamond-attains}. We do this over two steps:\\

\textbf{Step 1.} We first find approximations of $u$ and $\eta_l$ that allow us to take the limit as $l\to\infty$ on the sequence $\Prob_l$ and use weak convergence. Let $(c_m)_{m\geq 1}$ be a sequence of positive numbers. We note that we can (by approximating stopping times from the right and $U$ by a finite set) find a discrete approximation
\begin{align*}
  \hat u:=\sum_{i=1}^{M_m} \ett_{A^m_i}(\omega) \vecv_i,
\end{align*}
where $M_m\leq C/c_m^{k(d+1)}$ (recall that $d\geq 1$ is the dimension of $U$) and $\vecv_i\in D^k$, such that $A^m_i\in\mcF_{\tau_k}$ and $|\hat u-u|\leq c_m$ for all $\omega\in\Omega$.

Then, with $\eta^i_l:=\{t\geq \eta^i_{l-1}: S^{\vecv_i}_t-Y^{\vecv_i}_t\leq 1/l\}\wedge T$, we can, since $S^{\vecv_i}$ and $Y^{\vecv_i}$ are both uniformly continuous, repeat the argument in Step 1 in the proof of Theorem 3.3 in \cite{EkrenStopping} to find that there are continuous $[0,T]$-valued random variables $(\theta^i_l)_{l\geq 1}$ and $\mcF_T$-measurable sets $\Omega^i_l\subset\Omega$ such that
\begin{align*}
  \sup_{\Prob'\in\mcP(u)}\Prob'[(\Omega^i_l)^c]\leq c_l^{k(d+1)}2^{-l}\quad \text{and}\quad \eta^i_{l-1}-2^{-l}\leq \theta^i_l \leq \eta^i_{l+1}+2^{-l}\, \text{on}\,\Omega^i_l
\end{align*}
for $i=1,\ldots,M_l$. For $l\geq 1$, we define $\hat\eta_l:=\inf\{t\geq\hat \eta_{l-1}:S^{\hat u}_t-Y^{\hat u}_t\leq 1/l\}\wedge T$ with $\hat\eta_0:=\tau_k$ and note that $\hat\eta_l=\sum_{i=1}^{M_l} \ett_{A^l_i}(\omega) \eta^i_l$. Then $\Omega^{-}_l:=\{\omega\in\Omega:\hat\eta_l\leq \eta_{l-1}\}$ satisfies
\begin{align*}
  \Omega^-_l\subset \{\omega\in\Omega: S^{u}_{\hat\eta_l}(\omega)-Y^{u}_{\hat\eta_l}(\omega)\geq 1/(l-1)\}.
\end{align*}
Hence, as $S^{\hat u}_{\hat\eta_l}-Y^{\hat u}_{\hat\eta_l}=1/l$ on $[\hat\eta_l<T]$ by continuity, we have that
\begin{align*}
  \Omega^-_l\subset \{\omega\in\Omega: |Y^{u}_{\hat\eta_l}(\omega)-Y^{\hat u}_{\hat\eta_l}(\omega)|\geq 1/2l(l-1)\}\cup \{\omega\in\Omega: |S^{u}_{\hat\eta_l}(\omega)-S^{\hat u}_{\hat\eta_l}(\omega)|\geq 1/2l(l-1)\}
\end{align*}
However, there is a modulus of continuity function $\rho'$ such that
\begin{align*}
  \mcE^{u}[|S^{u}_{\hat\eta_l}(\omega)-S^{\hat u}_{\hat\eta_l}(\omega)|]&\leq \mcE^{u}\Big[(\rho_{Y,k}+\rho_{c,k})\Big(\sum_{j=1}^k\|\omega_{\tau_j+\cdot}-\omega_{\tau_j}\|_{c_l}+c_l\Big)\Big]\leq \rho'(c_l)
\end{align*}
and thus also
\begin{align*}
  \mcE^{u}[|Y^{u}_{\hat\eta_l}(\omega)-Y^{\hat u}_{\hat\eta_l}(\omega)|]&\leq \mcE^{u}\Big[\rho_{Y,k}\Big(\sum_{j=1}^k\|\omega_{\tau_j+\cdot}-\omega_{\tau_j}\|_{c_l}+c_l\Big)\Big]\leq \rho'(c_l)
\end{align*}
so that
\begin{align*}
  \sup_{\Prob\in\mcP(u)}\Prob[\Omega^-_l]\leq 4l(l-1)\rho'(c_l).
\end{align*}
Similarly, for $\Omega^{+}_l:=\{\omega\in\Omega:\hat\eta_l\geq \eta_{l+1}\}$ we have
\begin{align*}
  \sup_{\Prob\in\mcP(u)}\Prob[\Omega^+_l]\leq 4l(l+1)\rho'(c_l).
\end{align*}
We can thus choose $c_l$ such that $8l^2\rho'(c_l)\leq 2^{-l}$ for all $l\geq 1$ and letting $\hat\theta_l:=\sum_{i=1}^{M_l} \ett_{A^l_i}(\omega) \theta^i_l$ and $\hat\Omega_l:=\big(\cap_{i=1}^{M_l}\Omega^i_l\big)\cap (\Omega^{-}_l)^c\cap (\Omega^{+}_l)^c$ we find that
\begin{align*}
  \sup_{\Prob\in\mcP(u)}\Prob[(\hat\Omega_l)^c]\leq C2^{-l}\quad \text{and}\quad \eta_{l-2}-2^{-l}\leq \hat\theta_l \leq \eta_{l+2}+2^{-l}\, \text{on}\,\hat\Omega_l.
\end{align*}
Hence, since $\eta_{l}(\omega)\to\tau^\diamond(\omega)$ for all $\omega\in\Omega$, we get that $(\hat\theta_l)_{l\geq 1}$ is a sequence of random variables such that $\hat\theta_l$ is continuous on $A^l_i$ for $i=1,\ldots,M_l$ and by the Borel-Cantelli lemma, $\hat\theta_l\to\tau^\diamond$, $\Prob'$-a.s. for all $\Prob'\in\mcP(u)$.

Next, given $(c'_l)_{l\geq 1}$ there are $\Prob$-continuity sets $((D^l_i)_{i=1}^{M_l})_{l\geq 1}$ (\ie $\Prob[\partial D^l_i]=0)$ such that $D^l_i\in\mcF_{\tau_{k}}$ and $\Prob[A^l_i\Delta D^l_i]\leq c'_l$. Then, for each $\Prob'\in\mcP^{\Prob}_{\tau_k}(u)$, the sets in the double sequence $((D^l_i)_{i=1}^{M_l})_{l\geq 1}$ are also $\Prob'$-continuity sets and $\Prob'[A^l_i\Delta D^l_i]=\Prob[A^l_i\Delta D^l_i]\leq c'_l$. Moreover, $D^l_0:=(D^l_{1}\cup\cdots\cup D^l_{M_l})^c\in\mcF_{\tau_k}$ is also a $\Prob'$-continuity set for each $\Prob'\in\mcP^{\Prob}_{\tau_k}(u)$ and we can define
\begin{align*}
  u'_l:=\sum_{i=1}^{M_l} \ett_{D^l_i}(\omega) \vecv_i\quad\text{and}\quad \theta'_l:=\sum_{i=1}^{M_l} \ett_{D^l_i}(\omega) \theta^i_l.
\end{align*}
With $c'_l:=c_l^{k(d+1)}2^{-l}$ we thus find that $\theta'_l\to\tau^\diamond$, $\Prob'$-a.s.~for each $\Prob'\in\mcP^{\Prob}_{\tau_k}(u)$, moreover, $\theta'_l$ (and therefore also $Y^{u'_l}_{\theta'_l})$ is continuous on $D^l_i$ for $i=0,\ldots,M_l$.\\

\textbf{Step 2.} Using the approximation constructed in Step 1, we now show that \eqref{ekv:P-diamond-attains} holds. Since for any $t\in [\tau_k,T]$, we have
\begin{align}\label{ekv:up-approx}
  |Y^{u}_{t}-Y^{u'_l}_{t\vee \tau'_{l,k}}|\leq \rho_{Y,k}\Big(2\sum_{j=1}^k\|\omega_{\tau_j+\cdot}-\omega_{\tau_j}\|_{c_l}+c_l\Big)+2C_0\ett_{\Omega^\Delta_l},
\end{align}
where $\Omega^\Delta_l:=\cup_{i=1}^{M_l}A^l_i\Delta D^l_i$ and $\tau'_{l,k}$ is the time of the $k^\text{th}$ intervention in $u'_l$, with $\tau'_k:=T$ on $D^l_0$, we find that there is a $C>0$ such that for all $l\geq 2$, we have
\begin{align*}
  \E^{\Prob^\diamond}\big[Y^{u}_{\tau^\diamond}\big]&= \lim_{l\to\infty}\E^{\Prob^\diamond}\big[Y^{u'_l}_{\theta'_l}\big]=\lim_{l\to\infty}\lim_{m\to\infty}\E^{\Prob_m}\big[Y^{u'_l}_{\theta'_l}\big],
\end{align*}
where we have used the fact that, as $\Omega^\Delta_l\in\mcF_{\tau_k}$, we have
\begin{align}\label{ekv:delta-stab}
  \limsup_{l\to\infty}\limsup_{m\to\infty}\Prob_m[\Omega^\Delta_l]=\limsup_{l\to\infty}\Prob[\Omega^\Delta_l]=0.
\end{align}

On the other hand, by Theorem~\ref{thm:NutzZhang}.\emph{i)}, $Y^{u}_{\cdot\wedge\tau^\diamond}$ is a $\Prob'$-supermartingale for each $\Prob'\in\mcP(u)$ and we have that $\E^{\Prob_m}\big[Y^{u}_{\eta_l}\big]\geq \E^{\Prob_m}\big[Y^{u}_{\eta_m}\big]$ whenever $l\leq m$ implying that
\begin{align*}
  \liminf_{l\to\infty}\liminf_{m\to\infty}\E^{\Prob_m}\big[Y^{u}_{\eta_l}\big]\geq \liminf_{m\to\infty}\E^{\Prob_m}\big[Y^{u}_{\eta_m}\big]\geq \E^{\Prob}\big[Y^{u}_{\tau_k}\big],
\end{align*}
where, to reach the last inequality, we have used \eqref{ekv:prob_l-def} and the definition of $\eta_m$ to conclude that
\begin{align*}
  \E^{\Prob_m}\big[Y^u_{\eta_m}\big]\geq \E^{\Prob_m}\big[S^u_{\eta_m}\big]-1/m\geq \E^{\Prob}\big[\mcE^u_{\tau_k}\big[S^u_{\eta_m}\big]\big]-1/m-2^{-m}
\end{align*}
and the right hand side of the last inequality tends to $\E^{\Prob}\big[Y^u_{\tau_k}\big]$ as $m\to\infty$. Hence,
\begin{align*}
  \E^{\Prob}\big[Y^{u}_{\tau_k}\big]-\E^{\Prob^\diamond}\big[Y^{u}_{\tau^\diamond}\big]&\leq \liminf_{l\to\infty}\liminf_{m\to\infty}\E^{\Prob_m}\big[Y^{u}_{\eta_l}\big]- \lim_{l\to\infty}\lim_{m\to\infty}\E^{\Prob_m}\big[Y^{u}_{\theta'_l}\big]
  \\
  &\leq \limsup_{l\to\infty}\limsup_{m\to\infty}\E^{\Prob_m}\big[|Y^{u}_{\eta_l}-Y^{u'_l}_{\theta'_l}|\big]
  \\
  &\leq \limsup_{l\to\infty}\limsup_{m\to\infty}(\E^{\Prob_m}\big[|Y^{u'_l}_{\eta_l}-Y^{u'_l}_{\theta'_l}|\big]+\rho'(c_l) + 2C_0\Prob_m\big[\Omega^\Delta_l\big])
  \\
  &=\limsup_{l\to\infty}\limsup_{m\to\infty}\E^{\Prob_m}\big[|Y^{u'_l}_{\eta_l}-Y^{u'_l}_{\theta'_l}|\big]
\end{align*}
where we again use \eqref{ekv:up-approx} and \eqref{ekv:delta-stab}. Now, since $\eta_l$ is not continuous we cannot immediately proceed to take the limit. However, as
\begin{align*}
  \theta'_{l-2}-2^{2-l}\leq \eta_l\leq \theta'_{l+2}+2^{-2-l} \quad\text{on}\quad \Omega'_{l-2}\cap \Omega'_{l+2},
\end{align*}
where $\Omega'_l:=\hat\Omega_l\setminus \Omega^\Delta_l$, setting
\begin{align*}
  \psi_l:=\ett_{\cup_{i=1}^{M_l}D^l_i}\sup\big\{|Y^{u'_l}_{t}-Y^{u'_l}_{\theta'_l}|:t\in [\theta'_{l-2}-2^{2-l}, \theta'_{l+2}+2^{-2-l}]\cap[0,T]\big\},
\end{align*}
we find that $\psi_l$ is continuous on $D^l_i$ for $i=0,\ldots,M_l$ and thus
\begin{align*}
  \limsup_{l\to\infty}\limsup_{m\to\infty}\E^{\Prob_m}\big[\psi_l\big]=\limsup_{l\to\infty}\E^{\Prob^\diamond}\big[\psi_l\big]=0
\end{align*}
by dominated convergence under $\Prob^\diamond$. As
\begin{align*}
  \E^{\Prob_m}\big[|Y^{u'_l}_{\eta_l}-Y^{u'_l}_{\theta'_l}|\big]\leq \E^{\Prob_m}\big[\psi_l\big]+2C_0\Prob_m\big[(\Omega'_{l-2}\cap \Omega'_{l+2})^c\big]
\end{align*}
we thus find that
\begin{align*}
  \E^{\Prob}\big[Y^{u}_{\tau_k}\big]-\E^{\Prob^\diamond}\big[Y^{u}_{\tau^\diamond}\big]&\leq 2C_0\limsup_{l\to\infty}\limsup_{m\to\infty}\Prob_m\big[(\Omega'_{l-2}\cap \Omega'_{l+2})^c\big]
  \\
  &\leq C\limsup_{l\to\infty}(2^{-l}+\Prob\big[\Omega^\Delta_l\big])
  \\
  &=0
\end{align*}
from which the assertion follows since $Y^{u}_{\cdot\wedge\tau^\diamond}$ is a $\Prob^\diamond$-supermartingale.\qed\\

Since stopping before $\tau^\diamond$ in the above lemma is never optimal, this means that the $\Prob^\diamond$ in Lemma~\ref{lem:has-P-kernel-1} is optimal until $\tau^\diamond$. It remains to decide a continuation after $\tau^\diamond$ such that, at time $\tau^\diamond$, stopping is optimal.

\begin{lem}\label{lem:has-P-kernel-after}
Let $\mcP(t,\omega,\vecv)$ be weakly compact for all $(t,\omega,\vecv)\in\Lambda\times\mcD$ and let $u\in\mcU^k$. There is a $\mcF^*_{\tau_k}$-measurable kernel $\nu:\Omega\to\PrM(\Omega)$ such that $\nu(\omega)\in\mcP(\tau_k,\omega,u)$ and
\begin{align}\nonumber
&\inf_{\tau\in\mcT^{\tau_k(\omega)}}\E_{\tau_k(\omega)}^{\nu(\omega)}\Big[\ett_{[\tau=T]}\varphi(u) + \ett_{[\tau<T]}\inf_{b\in U}\big\{Y^{u\circ(\tau,b)}_\tau+c(u\circ(\tau,b))\big\}\Big](\omega)
\\
&=\inf_{\tau\in\mcT^{\tau_k(\omega)}}\mcE_{\tau_k(\omega)}^{u}\Big[\ett_{[\tau=T]}\varphi(u) + \ett_{[\tau<T]}\inf_{b\in U}\big\{Y^{u\circ(\tau,b)}_\tau+c(u\circ(\tau,b))\big\}\Big](\omega)
\end{align}
for all $\omega\in\Omega$.
\end{lem}

\noindent \emph{Proof.} Although the statement in the lemma differs from that of Lemma 4.16 in \cite{NutzZhang15}, the proof is almost identical and we give the main steps for the sake of completeness. We simplify notation by letting
\begin{align*}
  V(\vecv,\omega,\Prob):=&\inf_{\tau\in\mcT^{t_k}}\E_{t_k}^{\Prob}\Big[\ett_{[\tau=T]}\varphi(\vecv) + \ett_{[\tau<T]}\inf_{b\in U}\{Y^{\vecv\circ(\tau,b)}_\tau+c(\vecv\circ(\tau,b))\}\Big](\omega)
  \\
  =&\inf_{\tau\in\mcT}\E^{\Prob}\Big[\big(\ett_{[\tau\vee t_k=T]}\varphi(\vecv) + \ett_{[\tau\vee t_k<T]}\inf_{b\in U}\{Y^{\vecv\circ(\tau,b)}_\tau+c(\vecv\circ(\tau,b))\}\big)(\omega\otimes_{t_k}\cdot)\Big]
\end{align*}
for $\vecv\in D^k$. Then by Lemma 4.15 of \cite{NutzZhang15}, the map $\Prob\mapsto V(\vecv,\omega,\Prob)$ is upper semi-continuous for every $(\omega,\vecv)\in\Omega\times\mcD$. To use a measurable selection theorem we need to ascertain that the map $(\omega,\Prob)\mapsto V(u(\omega),\omega,\Prob)$ is Borel. Since $u\in\mcU^k$, we have that $\omega\mapsto u(\omega)$ is Borel and this will be obtained by showing that $(\vecv,\omega)\mapsto V(\vecv,\omega,\Prob)$ is Borel for any $\Prob\in\PrM(\Omega)$. However, for any $\tau\in \mcT$ we have that
\begin{align*}
  (\vecv,\omega)\mapsto\E^{\Prob}\Big[\big(\ett_{[\tau\vee t_k=T]}\varphi(\vecv) + \ett_{[\tau\vee t_k<T]}\inf_{b\in U}\{Y^{\vecv\circ(\tau,b)}_\tau+c(\vecv\circ(\tau,b))\}\big)(\omega\otimes_{t_k}\cdot)\Big]
\end{align*}
is Borel, by Fubini's theorem. Hence, taking the infimum over a countable set of stopping times that is dense outside of a $\Prob$-null set the assertion follows.\qed\\

\noindent \emph{Proof of \eqref{ekv:Pstar-is-opt}.} \textbf{Step 1.} We first construct a candidate for an optimal strategy. Note that by Theorem~\ref{thm:NutzZhang}, there is a $\Prob^*_0\in\mcP(\emptyset)$ such that
\begin{align}\nonumber
  Y_0&=\E^{\Prob^*_0}\Big[\ett_{[\tau^*_1=T]}\varphi(\emptyset) + \ett_{[\tau_1^*<T]}\{Y^{(\tau^*_1,\beta^*_1)}_{\tau^*_1}+c(\tau^*_1,\beta^*_1)\}\Big]
  \\
  &= \inf_{\tau\in\mcT}\E^{\Prob^*_0}\Big[\ett_{[\tau=T]}\varphi(\emptyset) + \ett_{[\tau<T]}\inf_{b\in U}\{Y^{(\tau,b)}_{\tau}+c(\tau,b)\}\Big].\label{ekv:dbl-first}
\end{align}
We fix $u\in\mcU$ and define $(\Prob^u_k)_{k\geq 1}$ iteratively where we start by setting $\Prob^u_0:=\Prob^*_0$. For $k=1,\ldots$, we construct $\Prob^{u}_{k}$ from $\Prob^u_{k-1}$ by first letting $\Prob^{u,\diamond}_k \in \mcP^{\Prob^u_{k-1}}_{\tau_{k-1}}([u]_{k-1})$ be as in the statement of Lemma~\ref{lem:has-P-kernel-1} with $\Prob=\Prob^u_{k-1}$ and $u\leftarrow [u]_k$, then there is by Lemma~\ref{lem:has-P-kernel-after} a $\mcF^*_{\tau^\diamond}$-measurable kernel $\nu^{u,k}:\Omega\to\PrM(\Omega)$ such that $\nu^{u,k}(\omega)\in\mcP(\tau_k,\omega,[u]_k)$ and
\begin{align}\nonumber
&\inf_{\tau\in\mcT^{\tau_k(\omega)}}\E_{\tau_k(\omega)}^{\nu^{u,k}(\omega)}\Big[\ett_{[\tau=T]}\varphi([u]_{k}) + \ett_{[\tau<T]}\inf_{b\in U}\{Y^{[u]_{k}\circ(\tau,b)}_\tau+c([u]_{k}\circ(\tau,b))\}\Big](\omega)
\\
&=\inf_{\tau\in\mcT^{\tau_k(\omega)}}\mcE_{\tau_k(\omega)}^{[u]_{k}}\Big[\ett_{[\tau=T]}\varphi([u]_{k}) + \ett_{[\tau<T]}\inf_{b\in U}\{Y^{[u]_{k}\circ(\tau,b)}_\tau+c([u]_{k}\circ(\tau,b))\}\Big](\omega)\label{ekv:verthm-P-att}
\end{align}
for all $\omega\in\Omega$ and $k\geq 1$. We can thus define the measure
\begin{align*}
    \Prob^u_{k}[A]:=\int\!\!\!\int (\ett_A)^{\tau_{k},\tilde \omega}(\omega')\hat\nu^{u,k}(d\omega';\tilde\omega)\Prob^{u,\diamond}_{k}[d\tilde\omega]
\end{align*}
for all $A\in\mcF$, where $\hat\nu^{u,k}$ is an $\mcF_{\tau_k}$-measurable kernel such that \eqref{ekv:verthm-P-att} holds $\Prob_k^{u,\diamond}$-a.s. Combining Lemma~\ref{lem:has-P-kernel-1} and Lemma~\ref{lem:has-P-kernel-after} we realize that our above construction extends \eqref{ekv:dbl-first} yielding that for each $k\geq 0$, we have
\begin{align}\nonumber
  \E^{\Prob^u_{k-1}}\big[Y^{[u]_k}_{\tau_k}\big]&=\E^{\Prob^u_{k}}\Big[\ett_{[\tau^\diamond=T]}\varphi([u]_k) + \ett_{[\tau^\diamond<T]}\inf_{b\in U}\{Y^{[u]_k\circ(\tau^\diamond,b)}_{\tau^\diamond}+c([u]_k\circ(\tau^\diamond,b))\}\Big]
  \\
  &=\inf_{\tau\in\mcT_{\tau_k}}\E^{\Prob^u_{k}}\Big[\ett_{[\tau=T]}\varphi([u]_k) + \ett_{[\tau<T]}\inf_{b\in U}\{Y^{[u]_k\circ(\tau,b)}_{\tau}+c(u\circ(\tau,b))\}\Big],\label{ekv:dbl-attainment}
\end{align}
where $\tau^\diamond:=\inf \big\{s \geq \tau_{k}:\:Y_s^{[u]_k}=\inf_{b\in U} \{Y^{[u]_k\circ(s,b)}_s+c([u]_k\circ(s,b))\}\big\}\wedge T$ and $\Prob^u_{-1}=\Prob_0$. Our candidate strategy is then $\Prob^{*,S}(u):=\Prob^u$, where $\Prob^u$ is such that (a subsequence of) $(\Prob^u_k)_{k\geq 0}$ converges weakly to $\Prob^u$. Since $\mcP^{\Prob^u_k}_{\tau_k}([u]_k)$ is weakly compact, we note that $\Prob^u=\Prob^u_k$ on $\mcF_{\tau_k}$ for each $k\geq 0$ and as $\Prob^u_k\in\mcP([u]_k)$ it follows by Assumption~\ref{ass:mcP} that $\Prob^u\in\mcP(u)$. Furthermore, it is clear by the above construction that the strategy $\Prob^{*,S}$ is non-anticipative, and we conclude that $\Prob^{*,S}\in\mcP^S$.\\

\textbf{Step 2.} Next, we show that if a subsequence of $(\Prob^{u^*}_k)_{k\geq 0}$ converges weakly to some $\Prob^*$, then $\Prob^*$ is an optimal response to the optimal impulse control $u^*$. We thus assume that $u$ in the previous step equals $u^*$. By \eqref{ekv:dbl-first} and the fact that $Y^{(T,\beta^*_1)}_{T}=\varphi(\emptyset)$ (see Remark~\ref{rem:@end}) we get that
\begin{align*}
  Y_0&=\E^{\Prob^*_0}\Big[\ett_{[\tau^*_1=T]}\varphi(\emptyset) + \ett_{[\tau_1^*<T]}\{Y^{(\tau^*_1,\beta^*_1)}_{\tau^*_1}+c(\tau^*_1,\beta^*_1)\}\Big]
  \\
  &= \E^{\Prob^*_0}\Big[ Y^{(\tau^*_1,\beta^*_1)}_{\tau^*_1}+\ett_{[\tau_1^*<T]}c(\tau^*_1,\beta^*_1)\Big].
\end{align*}
However, \eqref{ekv:dbl-attainment} gives that
\begin{align*}
  \E^{\Prob^*_0}\big[ Y^{[u^*]_1}_{\tau^*_1}\big]= \E^{\Prob^*_1}\Big[\ett_{[\tau^*_2=T]}\varphi(u^*) + \ett_{[\tau^*_2<T]}\{Y^{[u^*]_2}_{\tau^*_2}+c([u^*]_2)\}\Big].
\end{align*}
Repeating this procedure $k\geq 1$ times we find that
\begin{align*}
  Y^{\emptyset}_0&=\E^{\Prob^*_k}\big[\ett_{[\tau_{k+1}^*=T]}\varphi(u^*) + \sum_{j=1}^{N^*\wedge k+1} c([u^*]_j)+\ett_{[\tau^*_{k+1}<T]}Y^{[u^*]_{k+1}}_{\tau^*_{k+1}}\big].
\end{align*}
In fact, as $\Prob^*_k=\Prob^*_l$ on $\mcF_{\tau_{k+1}\wedge \tau_{l+1}}$ this holds when replacing the measure $\Prob^*_k$ with $\Prob^*_l$ whenever $l\geq k$ and we have that
\begin{align*}
  Y^{\emptyset}_0&=\lim_{k\to\infty}\lim_{l\to\infty}\E^{\Prob^*_l}\big[\ett_{[\tau_{k+1}^*=T]}\varphi(u^*) + \sum_{j=1}^{N^*\wedge k+1} c([u^*]_j)+\ett_{[\tau^*_{k+1}<T]}Y^{[u^*]_{k+1}}_{\tau^*_{k+1}}\big].
\end{align*}
By weak compactness of $\mcP^{\Prob^*_l}_{\tau^*_l}(u^*)$ for each $l\geq 0$, by possibly going to a subsequence, $\Prob^*_l$ converges weekly to some $\Prob^*\in \mcP(u^*)$ such that $\Prob^*=\Prob^*_l$ on $\mcF_{\tau^*_l}$ for each $l\geq 0$ and we find that
\begin{align*}
  Y^{\emptyset}_0&=\lim_{k\to\infty}\E^{\Prob^*}\big[\ett_{[\tau_{k+1}^*=T]}\varphi(u^*) + \sum_{j=1}^{N^*\wedge k+1} c([u^*]_j)+\ett_{[\tau^*_{k+1}<T]}Y^{[u^*]_{k+1}}_{\tau^*_{k+1}}\big]
  \\
  &=\E^{\Prob^*}\big[\varphi(u^*) + \sum_{j=1}^{N^*} c([u^*]_j)\big]
\end{align*}
by dominated convergence under $\Prob^*$.\\

\textbf{Step 3.} It remains to show that $u^*$ is an optimal impulse control under the strategy $\Prob^{*,S}$. For arbitrary $u\in\mcU$ we have
\begin{align*}
  Y^\emptyset_0&\leq\E^{\Prob^*_0}\Big[\ett_{[\tau_1=T]}\varphi(\emptyset) + \ett_{[\tau_1<T]}\{Y^{(\tau_1,\beta_1)}_{\tau_1}+c(\tau_1,\beta_1)\}\Big]
  \\
  &=\E^{\Prob^*_0}\Big[ Y^{(\tau_1,\beta_1)}_{\tau_1} + \ett_{[\tau_1<T]}c(\tau_1,\beta_1)\Big],
\end{align*}
by again using that $Y^{(T,\beta_1)}_{T}=\varphi(\emptyset)$. Now, \eqref{ekv:dbl-attainment} gives that
\begin{align*}
  \E^{\Prob^*_0}\big[ Y^{(\tau_1,\beta_1)}_{\tau_1}\big]\leq \E^{\Prob^u_1}\Big[\ett_{[\tau_2=T]}\varphi(u) + \ett_{[\tau_2<T]}\{Y^{[u]_2}_{\tau_2}+c([u]_2)\}\Big].
\end{align*}
Repeating $k$ times we find that
\begin{align*}
  Y^{\emptyset}_0&\leq\E^{\Prob^u_k}\big[\ett_{[\tau_{k+1}=T]}\varphi(u) + \sum_{j=1}^{N\wedge k+1} c([u]_j)+\ett_{[\tau_{k+1}<T]}Y^{[u]_{k+1}}_{\tau_{k+1}}\big].
\end{align*}
Arguing as in Step 2 gives that there is a subsequence under which $\Prob^u_k$ converges weakly to some $\Prob^u\in \mcP(u)$ such that $\Prob^u=\Prob^u_l$ on $\mcF_{\tau_l}$ for each $l\geq 0$ and then
\begin{align*}
  Y^{\emptyset}_0&\leq\liminf_{k\to\infty}\E^{\Prob^u}\Big[\ett_{[\tau_{k+1}=T]}\varphi(u) + \sum_{j=1}^{N\wedge k+1} c([u]_j)+\ett_{[\tau_{k+1}<T]}Y^{[u]_{k+1}}_{\tau_{k+1}}\Big].
\end{align*}
Now either $u$ makes infinitely many interventions on some set of positive size under $\Prob^u$ in which case the right hand side equals $+\infty$ or we can use dominated and monotone convergence under $\Prob^u$ to find that
\begin{align*}
  Y^{\emptyset}_0&\leq\E^{\Prob^u}\Big[\varphi(u) + \sum_{j=1}^{N} c([u]_j)\Big].
\end{align*}
Combined, this proves \eqref{ekv:Pstar-is-opt}.\qed\\

\begin{rem}
Theorem~\ref{thm:optim-ctrls} presumes existence of a $\mcF_{\tau_j}$-measurable minimizer $\beta^*$, we note that such a minimizer always exists if $U$ is a finite set. The canonical example under which the compactness assumption holds is when the uncertainty/ambiguity stems from the set of all It\^o processes $\int_0^t b_s+\int_0^t\sigma_sdB_s$ with $(b,\sigma\sigma^\top)\in A$ for some compact, convex set $A\subset \R^d\times \bbS_+$ (here $\bbS_+$ denotes the set of positive semi-definite symmetric $d\times d$-matrices).
\end{rem}


\section{Application to path-dependent zero-sum stochastic differential games\label{sec:sdg}}
In the present section we extend the literature on stochastic differential games by considering the application of the above developments to zero-sum games of impulse versus continuous control in a path-dependent setting. For simplicity we only consider the driftless setting as this is sufficient to capture the main features when applying the above results to controlled path-dependent SDEs.

Throughout this section, we only consider impulse controls for which $\tau_j(\omega)\to T$ as $j\to\infty$ for all $\omega\in\Omega$. We will make frequent use of the following decomposition of \cadlag paths:

\begin{defn}
For $x\in\bbD$ (the set of $\R^d$-valued \cadlag paths) we let $\mcJ(x)=(\eta_j,\Delta x_j)_{j=1}^{n(T)}$ where $\eta_j$ is the time of the $j^\text{th}$ jump of $x$, $\Delta x_j:=x_{\eta_j}-x_{\eta_j-}$ the corresponding jump (with $x_{0-}:=0$) and for each $t\in[0,t]$, $n(t)$ is the number of jumps of $x$ in $[0,t]$. Moreover, we let $\mcC(x)$ be the path without jumps, \ie $(\mcC(x))_s=x_s-\sum_{i=1}^{n(s)}\Delta x_i$. 
\end{defn}

\subsection{Problem formulation\label{subsec:pf}}
We let $\mcA$ be the set of all progressively measurable \cadlag processes $\alpha:=(\alpha_s)_{0\leq s\leq T}$, taking values in a bounded Borel-measurable subset $A$ of $\R^d$ and let $\mcA^S:\mcU\to \mcA$ be the corresponding set of non-anticipative strategies. We then consider the problem of showing that
\begin{align}
  \inf_{u\in U}\sup_{\alpha\in\mcA}J(u,\alpha)=\sup_{\alpha\in\mcA^S}\inf_{u\in U}J(u,\alpha)\label{ekv:pd-sdg-form}
\end{align}
where $J$ is the cost functional (recall that $\E$ is expectation with respect to $\Prob_0$, the probability measure under which $B$ is a Brownian motion)
\begin{align}\label{ekv:J-def}
J(u,\alpha):=\E\Big[\int_0^T\phi(t,X^{u,\alpha}_t)dt+\psi(X^{u,\alpha}_T) + \sum_{j=1}^{N(\mcC(X^{u,\alpha}))}(\ell(\tau_j(\mcC(X^{u,\alpha})),X^{[u]_{j-1},\alpha}_{\tau_j(\mcC(X^{u,\alpha}))},\beta_j(\mcC(X^{u,\alpha})))\Big].
\end{align}
We assume that $X^{u,\alpha}:=\limsup_{j\to\infty}X^{[u]_{j},\alpha}$, where $\limsup$ is taken componentwise, and $X^{[u]_{j}}$ is defined recursively as the $\Prob_0$-a.s.~unique solution to
\begin{align}
X^{[u]_{j},\alpha}_{t}&=\int_{0}^t\sigma(s,X^{[u]_{j},\alpha},\alpha_s)dB_s+\sum_{i=1}^{j\wedge N(\mcC(X^{u,\alpha}))}\ett_{[\tau_i(\mcC( X^{u,\alpha}))\leq t]}\Gamma(\tau_i(\mcC(X^{u,\alpha})),X^{[u]_{i-1},\alpha},\beta_i(\mcC(X^{u,\alpha}))),\label{ekv:forward-sde}
\end{align}
where $\sigma:[0,T]\times\bbD\times A\to\bbS_{++}$ (here $\bbS_{++}$ denotes the set of positive definite symmetric $d\times d$-matrices) and $\Gamma:[0,T]\times\bbD\times U\to\R^d$ satisfy the conditions in Assumption~\ref{ass:onSFDE} below. Examining \eqref{ekv:J-def} and \eqref{ekv:forward-sde}, we note that the control $u$ is implemented as a feedback control. In this framework, the intervention times are generally referred to as stopping rules rather than stopping times, following~\cite{KarSud01}, and the corresponding $u$ is referred to as a strategy (as in for example \cite{Sirbu14}) rather than a control as it depends on the history of $\alpha$ through the state process $X^{u,\alpha}$.

We will make frequent use of the following operator that acts as an inverse to $\mcC$:
\begin{defn}
For each $u\in\mcU$, we introduce the impulse operator $\mcI^u:\Omega\to\bbD$ defined as $\mcI^u(\omega):=\lim\sup_{k\to\infty}\mcI^{[u]_k}(\omega)$, where
\begin{align*}
(\mcI^{[u]_k}(\omega))_s :=\omega_s+\sum_{j=1}^{N(\omega)\wedge k}\ett_{[\tau_j(\omega)\leq s]}\Gamma(\tau_j(\omega),\mcI^{[u]_{j-1}}(\omega),\beta_j(\omega))
\end{align*}
for each $u\in\mcU$.
\end{defn}
We then let $\bar X^{u,\alpha}$ be the $\Prob_0$-a.s.~unique solution to the path-dependent SDE
\begin{align}\label{ekv:sde-barX}
\bar X^{u,\alpha}_t=\int_{0}^t\sigma(s,\mcI^u(\bar X^{u,\alpha}),\alpha_s)dB_s,
\end{align}
for all $(u,\alpha)\in\mcU\times\mcA$. The problem posed in \eqref{ekv:pd-sdg-form} can now be related to that of \eqref{ekv:sdg-form} by noting that
\begin{align}
J(u,\alpha)=\E^{\Prob^{u,\alpha}}\Big[\int_0^T\phi(t,(\mcI^{u}(\cdot))_t)dt+\psi((\mcI^{u}(\cdot))_T) + \sum_{j=1}^N\ell(\tau_j,(\mcI^{[u]_{j-1}}(\cdot)),\beta_j)\Big],
\end{align}
where $\Prob^{u,\alpha}:=\Prob_0\circ\big(\bar X^{u,\alpha}\big)^{-1}$. In the remainder of this section we first give suitable assumptions on the coefficients of the problem and then show that, under these assumptions, the problem stated here falls into the framework handled in the first part of the paper.

\subsection{Assumptions}

Throughout, we make the following assumptions on the coefficients in the above problem formulation:
\begin{ass}\label{ass:onSFDE}
For any $t,t'\in [0,T]$, $b,b'\in U$, $\xi \in\R^d$, $x,x'\in\bbD$ and $a\in A$ we have:
\begin{enumerate}[i)]
  \item\label{ass:onSFDE-Gamma} The function $\Gamma:[0,T]\times\bbD\times U\to\R^d$ is such that $(t,\omega)\mapsto \sigma (t,X,b)$ is progressively measurable whenever $X$ is progressively measurable. Moreover, there is a $C>0$ such that
  \begin{align*}
    |\Gamma(t,x,b)-\Gamma(t',x',b')|&\leq C(\textbf d_{n(t)}[(t,\mcC(x),\mcJ(x_{\cdot\wedge t})),(t',\mcC(x'),\mcJ(x'_{\cdot\wedge t'}))]+|b'-b|),
  \end{align*}
  whenever $n(t)=n'(t')$ (that is, the number of jumps of $x$ during $[0,t]$ and the number of jumps of $x'$ during $[0,t']$ agree).
  \item\label{ass:onSFDE-a-sigma} The coefficient  $\sigma:[0,T]\times\bbD\times A\to\bbS_{++}$ is such that $(t,\omega)\mapsto \sigma (t,X,\alpha)$ is progressively measurable (resp. \cadlagSTOP) whenever $X$ and $\alpha$ are progressively measurable (resp. \cadlagSTOP) and satisfies the growth condition
  \begin{align*}
    |\sigma(t,x,a)|&\leq C
  \end{align*}
  and the functional, resp. $L^2$, Lipschitz continuity
  \begin{align*}
    |\sigma(t,x',a)-\sigma(t,x,a)|&\leq C\sup_{s\leq t}|x'_s-x_s|,
    \\
    \int_0^t |\sigma(s,x',\alpha_s)-\sigma(s,x,\alpha_s)|^2ds&\leq C\int_{0}^t|x'_s-x_s|^2ds,
  \end{align*}
  for all $\alpha\in\mcA$.
  Moreover, for each $(t,x)\in[0,T]\times\bbD$, the map $a\mapsto\sigma(t,x,a)$ has a measurable inverse, \ie there is a $\sigma^\text{inv}:[0,T]\times\bbD\times\bbS_{++}\to A$, that satisfies the same measurability properties as $\sigma$, for which
  \begin{align*}
    \sigma^\text{inv}(t,x,\sigma(t,x,a))=a
  \end{align*}
  for all $(t,x,a)\in [0,T]\times \bbD \times A$.
  \item\label{ass:onSFDE-phi} The running reward $\phi:[0,T]\times \R^d\to\R$ is Borel-measurable, continuous in $a$, uniformly bounded and uniformly continuous in the second variable uniformly in the other variables.
  \item\label{ass:onSFDE-psi} The terminal reward $\psi:\R^d\to\R$ is uniformly bounded and uniformly continuous.
  \item\label{ass:onSFDE-c} The intervention cost $\ell:[0,T]\times \R^d\times U\to \R_+$ is uniformly bounded, uniformly continuous and satisfies
  \begin{align*}
    \ell(t,\xi,b)\geq\delta >0.
  \end{align*}
\end{enumerate}
\end{ass}

\subsection{Value of the game}
We now provide a solution to the problem stated in Section~\ref{subsec:pf} by showing that, under the above assumptions, it is a special case of the general impulse control problem under non-linear expectation treated in sections \ref{sec:prel}-\ref{sec:verif-thm}.

For $(t,\omega)\in\Lambda$, we extend the definition of $\bar X^{u,\alpha}$ by letting $\bar X^{u,\alpha,t,\omega}$ be the $\Prob_0$-a.s.~unique solution to the path-dependent SDE
\begin{align*}
\bar X^{u,\alpha,t,\omega}_s=\omega_t+\int_{0}^s\sigma(t+r,\mcI^u(\omega\otimes_t\bar X^{u,\alpha,t,\omega}),\alpha_{t+r})dB_r.
\end{align*}
We then let $\mcP(t,\omega,u):=\cup_{\alpha\in\mcA}\{\Prob_0\circ(\bar X^{u,\alpha,t,\omega})^{-1}\}$ and get, in particular, that $\mcP(u):=\cup_{\alpha\in \mcA}\{\Prob_0\circ (\bar X^{u,\alpha})^{-1}\}$. In this setting it is straightforward to show that $\mcP(t,\omega,u)$ corresponds to a random $G$-expectation \cite{Nutz-G-exp} and thus satisfies properties 1-3 of Assumption~\ref{ass:mcE-general} by the results of \cite{NeufeldNutz13}. Moreover, by non-anticipativity and continuity of $\bar X^{u,\alpha}$, the fourth property also follows and we conclude that Assumption~\ref{ass:mcE-general} holds. Note also that Assumption~\ref{ass:mcE-bound} holds since $\sigma$ is bounded.

It remains to show that Assumptions \ref{ass:costs} and \ref{ass:mcE-consist} hold. For this we introduce for each $u\in\mcU^{t}$, the controlled process $\hat X^{t,\omega,u,\alpha}:=\limsup_{j\to\infty}\hat X^{t,\omega,[u]_{j},\alpha}$, where
\begin{align}
\hat X^{t,\omega,[u]_{j},\alpha}_{\cdot}&=(\omega\otimes_t\int_{t}^\cdot\sigma(s,\hat X^{t,\omega,[u]_{j},\alpha},\alpha_s)dB_s)+\sum_{i=1}^{j\wedge N(\cdot)(0\otimes_t B)}\Gamma(\tau_i(0\otimes_t B),\hat X^{t,\omega,[u]_{i-1},\alpha},\beta_i(0\otimes_t B)),\label{ekv:hatX-sde}
\end{align}
has $u$ implemented as an open loop strategy. Moreover, we define the control set $\mcA_t$ as the set of all $(\alpha_s)_{t\leq s\leq T}$ with $\alpha\in \mcA$ and $(\alpha_s(\omega))_{t\leq s\leq T}$ independent of $\omega|_{[0,t)}$ and for each $\vecv=(t_j,b_j)_{1\leq j\leq k},\vecv'=(t'_j,b'_j)_{1\leq j\leq k}\in D^k$, we introduce the set $\Upsilon^{\vecv,\vecv'}:=\cup_{j=1}^k[t_j\wedge t'_j,t_j\vee t'_j)$. The core result we need is the following stability estimate:
\begin{prop}\label{prop:momentX}
Under Assumption~\ref{ass:onSFDE} there is for each $k,\kappa\geq 0$ a modulus of continuity $\rho_{X,\kappa+k}$ such that
\begin{align}
\E\Big[\sup_{s\in [0,T]\cap \Upsilon^{\vecv,\vecv'}}|\hat X^{t,\omega,\vecv\circ u,\alpha}_{s}-\hat X^{t,\omega',\vecv'\circ u,\alpha}_s|\Big]\leq \rho_{X,\kappa+k}(\mathbf d_\kappa[(t,\omega,\vecv),(t,\omega',\vecv')]),
\end{align}
for all $(t,\omega,\omega',\vecv,\vecv') \in [0,T]\times\Omega^2\times (D^\kappa)^2$ and $(u,\alpha)\in \mcU^{k}_t\times\mcA_t$.
\end{prop}

\noindent\emph{Proof.} To simplify notation we let $v:=\vecv\circ u=(\eta_j,\gamma_j)_{j=1}^{k+\kappa}$ and $v':=\vecv'\circ u=(\eta'_j,\gamma'_j)_{j=1}^{k+\kappa}$ and set $X^l:=\hat X^{t,\omega,[v]_l,\alpha}$ and $\tilde X^l:=\hat X^{t,\omega',[v']_l,\alpha}$ for $l=0,\ldots,\kappa+k$. Moreover, we let $\delta X^l:=X^l-\tilde X^l$ and set $\delta X:=\delta  X^{\kappa+k}$. For $s\in[\bar \eta_l,T]$ with $\bar \eta_l:=\eta_l\vee \eta'_l$ we have
\begin{align*}
\delta  X^{l}_{s}&=(\omega\otimes_t\int_t^\cdot\sigma(r,X^{l},\alpha_r)dB_r)_s- (\omega'\otimes_{t}\int_{t}^\cdot\sigma(r,\tilde X^{l},\alpha_r)dB_r)_s
\\
&\quad+\sum_{j=1}^l\Gamma(\eta_j,X^{j},\gamma_j)-\Gamma(\eta'_j,\tilde X^{j},\gamma'_j).
\end{align*}
Now, as $X^l=X^j$ on $[0,\eta_{j+1})$ for $j\leq l$, we have
\begin{align*}
  &|\Gamma(\eta_j,X^{j},\gamma_j)-\Gamma(\eta'_j,\tilde X^{j},\gamma'_j)|\leq C(\mathbf d_{j-1}[(\eta_j,\mcC(X^{j-1}),\mcJ( X^{j-1})),(\eta'_j,\mcC(\tilde X^{j-1}),\mcJ(\tilde X^{j-1}))]+|\gamma_j-\gamma'_j|)
  \\
  &\leq C(\sum_{i=1}^j\sup_{s\in [0,T]}|(\omega\otimes_t\int_t^\cdot\sigma(r,X^{l},\alpha_r)dB_r)_{s\wedge\eta_i}- (\omega'\otimes_{t}\int_{t}^\cdot\sigma(r,\tilde X^{l},\alpha_r)dB_r)_{s\wedge\eta'_i}|
  \\
  &\quad + \sum_{i=1}^{j-1}|\Gamma(\eta_i,X^{i},\gamma_i)-\Gamma(\eta'_i,\tilde X^{i},\gamma'_i)|+|\eta_j-\eta'_j|+|\gamma_j-\gamma'_j|)
\end{align*}
and induction gives that
\begin{align*}
  &|\Gamma(\eta_j,X^{j},\gamma_j)-\Gamma(\eta'_j,\tilde X^{j},\gamma'_j)|
  \\
  &\leq C(\sum_{i=1}^j\sup_{s\in [0,T]}|(\omega\otimes_t\int_t^\cdot\sigma(r,X^{l},\alpha_r)dB_r)_{s\wedge\eta_i}- (\omega'\otimes_{t}\int_{t}^\cdot\sigma(r,\tilde X^{l},\alpha_r)dB_r)_{s\wedge\eta'_i}| + |[v]_j-[v']_j|).
\end{align*}
Combined, we find that
\begin{align}\nonumber
|\delta  X^{l}_{s}|&\leq C(\sum_{i=1}^{l}\sup_{s'\in [0,s]}|(\omega\otimes_t\int_t^\cdot\sigma(r,X^{l},\alpha_r)dB_r)_{s'\wedge\eta_i}- (\omega'\otimes_{t}\int_{t}^\cdot\sigma(r,\tilde X^{l},\alpha_r)dB_r)_{s'\wedge\eta'_i}|
\\
&\quad + \sup_{s'\in [0,s]}|(\omega\otimes_t\int_t^\cdot\sigma(r,X^{l},\alpha_r)dB_r)_{s'} - (\omega'\otimes_{t}\int_{t}^\cdot\sigma(r,\tilde X^{l},\alpha_r)dB_r)_{s'}| + |[v]_l-[v']_l|).\label{ekv:comb}
\end{align}
For $l=0$ we see that
\begin{align*}
\sup_{s'\in [0,s]}|\delta  X^{0}_{s'}|&\leq \sup_{s'\in [0,s]}|(\omega\otimes_t\int_t^\cdot\sigma(r,X^{0},\alpha_r)dB_r)_{s'}-(\omega'\otimes_{t}\int_{t}^\cdot\sigma(r,\tilde X^{0},\alpha_r)dB_r)_{s'}|.
\end{align*}
Squaring and taking expectations on both sides and using the Burkholder-Davis-Gundy inequality then gives
\begin{align*}
\E\Big[\sup_{s'\in [0,s]}|\delta  X^{0}_{s'}|^2\Big]&\leq C(\|\omega-\omega'\|^2_t+\E\Big[\int_{t}^s|\sigma(r,X^{0},\alpha_r)-\sigma(r,\tilde X^{0},\alpha_r)|^2dr\Big])
\\
&\leq C(\|\omega-\omega'\|^2_t+\E\Big[\int_{t}^s|\delta X^{0}_r|^2dr\Big]),
\end{align*}
where we have used the $L^2$-Lipschitz condition and the uniform bound on $\sigma$ to arrive at the last inequality. Gr\"onwall's lemma now gives that
\begin{align*}
\E\Big[\sup_{s\in [0,T]}|\delta  X^{0}_{s}|^2\Big]&\leq C\textbf d[(t,\omega,\emptyset),(t,\omega',\emptyset)]^2
\end{align*}
We proceed by induction and assume that there is a $C$ such that
\begin{align*}
\E\Big[\sup_{s\in [0,T]\cap \Upsilon^{[\vecv]_l,[\vecv']_l}}|\delta  X^{l}_{s}|^2\Big]&\leq C(\textbf d[(t,\omega,[\vecv]_l),(t,\omega',[\vecv']_l)]^2+|[\vecv]_l-[\vecv']_l|).
\end{align*}
But then, repeating the above argument and using \eqref{ekv:comb}, gives that
\begin{align*}
\E\Big[\sup_{s\in [0,T]\cap \Upsilon^{[\vecv]_{l+1},[\vecv']_{l+1}}}|\delta  X^{l+1}_{s}|^2\Big]&\leq C(\E\Big[\sup_{s\in [0,T]\cap \Upsilon^{[v]_l,[v']_l}}|\delta  X^{l}_{s}|^2\Big]+\|\omega_{\cdot\wedge \eta_{l+1}}-\omega'_{\cdot\wedge \eta'_{l+1}}\|_T
\\
&\quad +|\eta_{l+1}-\eta'_{l+1}|^2+|\gamma_{l+1}-\gamma'_{l+1}|^2+|\eta_{l+1}-\eta'_{l+1}|)
\\
&\leq C(\textbf d[(t,\omega,[\vecv]_{l+1}),(t,\omega',[\vecv']_{l+1})]^2+|[\vecv]_{l+1}-[\vecv']_{l+1}|).
\end{align*}
Since $|v-v'|\leq k|\vecv-\vecv'|$ and $\Upsilon^{v,v'}=\Upsilon^{\vecv,\vecv'}$ we find by induction that
\begin{align*}
\E\Big[\sup_{s\in [0,T]\cap \Upsilon^{\vecv,\vecv'}}|\delta  X^{k+\kappa}_{s}|^2\Big]&\leq C(\textbf d[(t,\omega,\vecv),(t,\omega',\vecv')]^2+|\vecv-\vecv'|)
\end{align*}
and the statement of the proposition follows by Jensen's inequality.\qed\\

\begin{rem}
Note that the constant $C>0$ in the proof of Proposition~\ref{prop:momentX} generally depends on $k+\kappa$ and may tend to infinity as either $k$ or $\kappa$ goes to infinity. This is the main reason that we needed to go through the process of truncating the maximal number of interventions to obtain the dynamic programming principle in Section~\ref{sec:DynPP} and to prove that the game has a value in Section~\ref{sec:value}.
\end{rem}


To translate our problem to the general framework we define
\begin{align}\label{ekv:varphi-def}
\varphi(\omega,\vecv):=\int_0^T\phi(t,(\mcI^\vecv(\omega))_t)dt+\psi((\mcI^\vecv(\omega))_T)
\end{align}
and
\begin{align}\label{ekv:c-def}
c(\omega,(t_i,b_i)_{i=1}^j):=\ell(t_j,(\mcI^{(t_i,b_i)_{i=1}^{j-1}}(\omega))_{t_j},b_j).
\end{align}

\begin{lem}
The maps $\varphi$ and $c$ defined in \eqref{ekv:varphi-def}, resp.~\eqref{ekv:c-def} satisfy Assumption~\ref{ass:costs}.
\end{lem}

\noindent \emph{Proof.} Uniform boundedness of $\varphi$ (resp. $c$) is immediate from the boundedness of $\phi$ and $\psi$ and \eqref{ekv:varphi-def} (resp.~$\ell$ and \eqref{ekv:c-def}).

Letting $t=T$ in Proposition~\ref{prop:momentX} we find that, with $\vecv,\vecv'\in D^k$, we have
\begin{align*}
\sup_{s\in [0,T]\setminus\Upsilon^{\vecv,\vecv'}}|(\mcI^{\vecv}(\omega)-\mcI^{\vecv'}(\omega'))_s| & \leq \rho_{X,k}(\mathbf d[(T,\omega,\vecv),(T,\omega',\vecv')]).
\end{align*}
Since $\phi$ is bounded and the Lebesgue measure of the set $\Upsilon^{\vecv,\vecv'}$ is bounded by $|\vecv-\vecv'|$, we conclude that
\begin{align*}
|\varphi(\omega',\vecv')-\varphi(\omega,\vecv)|&\leq C|\vecv'-\vecv|+(T\rho_{\phi}+\rho_{\psi})\big(\rho_{X,k}(\mathbf d[(T,\omega,\vecv),(T,\omega',\vecv')])\big)
\\
&=:\rho_{\varphi,k}(\mathbf d[(T,\omega,\vecv),(T,\omega',\vecv')]),
\end{align*}
where $\rho_\phi$ and $\rho_\psi$ are bounded modulus of continuities for $\phi$ and $\psi$, respectively. Moreover, as
\begin{align*}
|(\mcI^{(t_i,b_i)_{i=1}^{k-1}}(\omega))_{t_k}-(\mcI^{(t'_i,b'_i)_{i=1}^{k-1}}(\omega'))_{t'_k}| & \leq \rho_{X,k}(\mathbf d[(t_k,\omega,\vecv),(t'_k,\omega',\vecv')]),
\end{align*}
we find that
\begin{align*}
|c(\omega',\vecv')-c(\omega,\vecv)|&\leq \rho_{\ell}(\rho_{X,k}(\mathbf d[(t_k,\omega,\vecv),(t'_k,\omega',\vecv')]))
\end{align*}
and we conclude that $c$ satisfies Assumption~\ref{ass:costs} with $\rho_{c,k}(\cdot)=\rho_\ell(\rho_{X,k}(\cdot))$, where $\rho_\ell$ is a bounded modulus of continuity for $\ell$.\qed\\

Finally, we have the following result:
\begin{lem}
Assumption~\ref{ass:mcE-consist} holds in this setting.
\end{lem}

\emph{Proof.} We fix $k\geq 0$, $u\in\mcU^{t,k}$ and $\omega'\in\Omega$ and note that $x\to \mcI^u(x)$ maps progressively measurable (resp. \cadlagSTOP) processes to progressively measurable (resp. \cadlagSTOP) processes and then so does the composition $\sigma^u:=\sigma(\cdot,\mcI^u(\cdot),\cdot)$.  We can thus repeated the argument in the appendix to \cite{NutzZhang15} to show that there is a $\mcF_t\otimes\mcF$-measurable map $(\omega,\tilde\omega)\mapsto u_\omega(\tilde\omega)$ such that $u_\omega\in\mcU^{t,k}$ for all $\omega\in\Omega$ and
\begin{align*}
  u_{\omega}(0\otimes_t\bar X^{u,\alpha,t,\omega})=u(0\otimes_t\bar X^{u,\alpha,t,\omega'}),\qquad \Prob_0\text{-a.s.}
\end{align*}
Now, since the set of open loop controls is at least as rich as the set of feedback controls, we have for $\vecv,\vecv'\in D^\kappa$, that
\begin{align*}
&\mcE^{\vecv'\circ u}_{t}\big[\varphi(\vecv'\circ u) + \sum_{j=1}^N c(\vecv'\circ [u]_j)\big](\omega')-\mcE^{\vecv\circ u_\omega}_t\big[\varphi(\vecv\circ u_\omega) + \sum_{j=1}^N c(\vecv\circ [u_\omega]_j)\big](\omega)
\\
&\leq \sup_{u\in\mcU^{t,k}}\sup_{\alpha\in\mcA_t}\E\Big[\int_0^T|\phi(s,\hat X^{t,\omega',\vecv'\circ u,\alpha}_s)-\phi(s,\hat X^{t,\omega,\vecv\circ u,\alpha}_s)|ds+|\psi(\hat X^{t,\omega',\vecv'\circ u,\alpha}_T)-\psi(\hat X^{t,\omega,\vecv\circ u,\alpha}_T)|
\\
&\quad +\sum_{j=1}^N|\ell(\tau_j\vee t'_\kappa,\hat X^{t,\omega',\vecv'\circ [u]_{j-1},\alpha}_{\tau_j\vee t'_\kappa},\beta_j)-\ell(\tau_j\vee t_\kappa,\hat X^{t,\omega,\vecv\circ [u]_{j-1},\alpha}_{\tau_j\vee t_\kappa},\beta_j)|\Big]
\end{align*}
Concerning the first two terms we have for any $u\in\mcU^{t,k}$, with $\rho'$ the concave hull of $\rho_\varphi+\rho_\psi$,
\begin{align*}
&\sup_{\alpha\in\mcA_t}\E\Big[\int_0^T|\phi(s,\hat X^{t,\omega',\vecv'\circ u,\alpha}_s)-\phi(s,\hat X^{t,\omega,\vecv\circ u,\alpha}_s)|ds+|\psi(\hat X^{t,\omega',\vecv'\circ u,\alpha}_T)-\psi(\hat X^{t,\omega,\vecv\circ u,\alpha}_T)|\Big]
\\
&\leq C(|\vecv-\vecv'|+\sup_{\alpha\in\mcA_t}\E\Big[\sup_{s\in [0,T]\cap \Upsilon^{\vecv,\vecv'}}\rho'(|\hat X^{t,\omega',\vecv'\circ u,\alpha}_s-\hat X^{t,\omega,\vecv\circ u,\alpha}_s|)\Big]
\\
&\leq C(|\vecv-\vecv'|+\rho'(\sup_{\alpha\in\mcA_t}\E\Big[\sup_{s\in [0,T]\cap \Upsilon^{\vecv,\vecv'}}|\hat X^{t,\omega',\vecv'\circ u,\alpha}_s-\hat X^{t,\omega,\vecv\circ u,\alpha}_s|\Big])
\end{align*}
Similarly, we have that
\begin{align*}
&\sup_{\alpha\in\mcA_t}\E\Big[\ell(\tau_j\vee t'_\kappa,\hat X^{t,\omega',\vecv'\circ [u]_{j-1},\alpha}_{\tau_j\vee t'_\kappa},\beta_j)-\ell(\tau_j\vee t_\kappa,\hat X^{t,\omega,\vecv\circ [u]_{j-1},\alpha}_{\tau_j\vee t_\kappa},\beta_j)|\Big]
\\
&\leq \sup_{\alpha\in\mcA_t}\E\Big[\rho_\ell(|t'_\kappa-t_\kappa|+|\hat X^{t,\omega',\vecv'\circ u,\alpha}_{t'_\kappa}-\hat X^{t,\omega,\vecv\circ u,\alpha}_{t_\kappa}|)\Big]
\\
&\leq \rho''(|t'_\kappa-t_\kappa|)+\rho''(\sup_{\alpha\in\mcA_t}\E\Big[|\hat X^{t,\omega',\vecv'\circ u,\alpha}_{t'_\kappa}-\hat X^{t,\omega,\vecv\circ u,\alpha}_{t_\kappa}|\Big])
\end{align*}
for some modulus of continuity $\rho''$. The result now follows by taking the supremum over all $u\in\mcU$ and using Proposition~\ref{prop:momentX} combined with a standard estimate for SDEs.\qed\\


\bibliographystyle{plain}
\bibliography{ImpNonlinExp_ref}
\end{document}